\documentclass[11pt]{article}%
\usepackage{amssymb}
\usepackage{amsmath}
\usepackage{amsfonts}
\usepackage{graphicx}
\usepackage[onehalfspacing]{setspace}
\usepackage{geometry}%
\providecommand{\U}[1]{\protect\rule{.1in}{.1in}}
\providecommand{\U}[1]{\protect\rule{.1in}{.1in}}
\newtheorem{theorem}{Theorem}

\newtheorem{lemma}[theorem]{Lemma}

\newtheorem{proposition}[theorem]{Proposition}
\newtheorem{remark}[theorem]{Remark}

\begin{document}

\title{Local Asymptotic Normality of the spectrum of high-dimensional spiked F-ratios}
\author{Prathapasinghe Dharmawansa, Iain M. Johnstone, and Alexei Onatski.}
\maketitle

\begin{abstract}
We consider two types of \textit{spiked} multivariate F distributions: a
scaled distribution with the scale matrix equal to a rank-one perturbation of
the identity, and a distribution with trivial scale, but rank-one
non-centrality. The norm of the rank-one matrix (\textit{spike}) parameterizes
the joint distribution of the eigenvalues of the corresponding F matrix. We
show that, for a \textit{spike} located above a \textit{phase transition}
threshold, the asymptotic behavior of the log ratio of the joint density of
the eigenvalues of the F matrix to their joint density under a local deviation
from this value depends only on the largest eigenvalue $\lambda_{1}$.
Furthermore, $\lambda_{1}$ is asymptotically normal, and the
\textit{statistical experiment} of observing all the eigenvalues of the F
matrix converges in the Le Cam sense to a \textit{Gaussian shift experiment}
that depends on the asymptotic mean and variance of $\lambda_{1}$. In
particular, the best statistical inference about a sufficiently large
\textit{spike} in the local asymptotic regime is based on the largest
eigenvalue only. As a by-product of our analysis, we establish joint
asymptotic normality of a few of the largest eigenvalues of the multi-spiked F
matrix when the corresponding spikes are above the \textit{phase transition}
threshold.\medskip\newline\textsc{Key words}: Spiked F-ratio, Local Asymptotic
Normality, multivariate F distribution, phase transition, super-critical
regime, asymptotic normality of eigenvalues, limits of statistical experiments.

\end{abstract}

\section{Introduction}

In this paper we establish the \textit{Local Asymptotic Normality}
(\textit{LAN}) of the \textit{statistical experiments} of observing the
eigenvalues of the F-ratio,\ $B^{-1}A,$ of two high-dimensional independent
\textit{Wishart} matrices, $A$ and $B$. We consider two situations. First,
both $A$ and $B$ are central Wisharts with dimensionality and degrees of
freedom that grow proportionally, and with the covariance parameters that
differ by a matrix of rank one. Second, $A$ and $B$ have the same covariance
parameter, but $A$ is a non-central Wishart with the non-centrality parameter
of rank one. In both cases, the joint distribution of the eigenvalues of
$B^{-1}A$ depends on the norm of the rank-one matrix, which we call a
\textit{spike}. We find that the considered \textit{statistical experiments}
are \textit{LAN} under a local parameterization of the \textit{spike} when the
locality is above a \textit{phase transition} threshold.

Many classical multivariate statistical tests are based on the eigenvalues of
F-ratio matrices. For example, all tests of the equality of two covariance
matrices and of the \textit{general linear hypothesis} in the
\textit{Multivariate Linear Model} described in Muirhead's (1982) chapters 8
and 10 are of this form. Contemporaneous statistical applications often
require the dimensionality of the F-ratio and its degrees of freedom be large
and comparable. Therefore, we consider the asymptotic regime where the
dimensionality and the degrees of freedom diverge to infinity at the same rate.

Our requirement that the parameters of the two Wisharts differ by a rank-one
matrix can be linked to situations where the alternative hypothesis is
characterized by the presence of one factor or signal, which is absent from
the data under the null. Inference conditional on factors requires considering
non-central F-ratios, whereas the unconditional inference leads to F-ratios
with unequal covariances.

The main result of this paper can be summarized as follows. We show that the
asymptotic behavior of the log ratio of the joint density of the eigenvalues
of $B^{-1}A,$ which corresponds to a sufficiently large value of the
\textit{spike}, to their joint density under a local deviation from this value
depends only on the largest eigenvalue $\lambda_{1}$. Furthermore,
$\lambda_{1}$ is asymptotically normal, and the \textit{statistical
experiment} of observing all the eigenvalues of $B^{-1}A$ converges in the Le
Cam sense to a \textit{Gaussian shift experiment} that depends on the
asymptotic mean and variance of $\lambda_{1}$. In particular, the best
statistical inference about a sufficiently large \textit{spike} in the local
asymptotic regime is based on the largest eigenvalue only.

We derive an explicit formula for the \textit{phase transition} threshold
demarcating the area of the sufficiently large \textit{spikes}. In a general
framework, where the parameters of $A$ and $B$ may differ by a matrix $\Delta$
of a finite rank, we show that, when the norm of $\Delta$ is below the
threshold, any finite number of the largest eigenvalues of $B^{-1}A$ almost
surely converge to the upper boundary of the support of the limiting spectral
distribution of $B^{-1}A,$ derived by Wachter (1980). In contrast, when $m$ of
the largest eigenvalues of $\Delta$ are above the threshold, we find that the
$m$ of the largest eigenvalues of $B^{-1}A$ almost surely converge to
locations strictly above the upper boundary of Wachter's distribution, and
that their local fluctuations about these limits are asymptotically jointly normal.

In a setting of two independent and not necessarily normal samples, the
\textit{phase transition} phenomenon has been studied in Nadakuditi and
Silverstein (2010). They obtain a formula for the threshold, and establish the
almost sure limits of the $m$ largest eigenvalues for the case where $\Delta$
describes the difference between covariance matrices of the two samples. The
limiting distribution of fluctuations above the threshold is described in
their paper as an open problem. Our paper solves this problem for the case of
two normal samples.

The \textit{phase transition} phenomenon for a single Wishart matrix has also
been a subject of active recent research. Baik et al (2005) study the joint
distributions of a few of the largest eigenvalues of \textit{complex Wisharts}
with spiked covariance parameters. They derive the asymptotic distributions of
a few of the largest eigenvalues, which turn out to be different depending on
whether the sizes of the corresponding spikes are below, at, or above a
\textit{phase transition} threshold, the situations often referred to as the
\textit{sub-critical}, \textit{critical}, and \textit{super-critical regimes}.

Similar transition takes place for real Wisharts. Paul (2007) establishes
asymptotic normality of the fluctuations of a few of the largest eigenvalues
in the \textit{super-critical regime} of the real case. F\'{e}ral and
P\'{e}ch\'{e} (2009), Benaych-Georges et al (2011) and Bao et al (2014a) show
that the fluctuations in the \textit{sub-critical} real case have the
\textit{Tracy-Widom distribution}, while Mo (2012) and Bloemendal and
Vir\`{a}g (2011, 2013) establish the asymptotic distribution of a different
type in the critical regime. In a setting of two normal samples, Bao et al
(2014b) study the almost sure limits of the \textit{sample canonical
correlations} when the \textit{population canonical correlations} are below
and when they are above a \textit{phase transition} threshold.

Our results on the joint asymptotic normality of the largest eigenvalues in
the \textit{super-critical regime} for F-ratios can be used to make
statistical inference about the eigenvalues of the \textquotedblleft
ratio\textquotedblright\ of the population covariances of $A$ and $B$, or the
eigenvalues of the non-centrality parameter of $A$. The estimates of these
eigenvalues play important role in MANOVA and the \textit{discriminant
analysis}, and can also be used in constructing modified model selection
criteria as discussed in Sheena et al (2004). Further, they may be important
in as diverse applications as constructing genetic selection indices and
describing a degree of financial turbulence (see Hayes and Hill (1981), and
Kritzman and Li (2010)).

We expect that our asymptotic normality results can be extended to the case of
the \textquotedblleft ratio\textquotedblright\ of two sample covariance
matrices constructed from non-normal samples. In the one-sample case, such an
extension of Paul's (2007) asymptotic normality results has been done in Bai
and Yao (2008). In this paper, we focus on normal data. This focus is dictated
by our main goal: establishing the \textit{LAN} property of the statistical
experiments of observing the eigenvalues of $B^{-1}A$. To reach this goal, we
derive an asymptotic approximation to a log likelihood process by representing
it in the form of a contour integral, and applying the \textit{Laplace
approximation} method. The explicit form of the joint distribution of the
eigenvalues of $B^{-1}A$ is known only in the normal case, and we need such an
explicit form for our analysis.

A decision-theoretic approach to the finite sample estimation of the
eigenvalues of the \textquotedblleft ratio\textquotedblright\ of the
population covariances of $A$ and $B$, or the eigenvalues of the
non-centrality parameter of $A$ was taken in many previous studies (see Sheena
et al (2004), Bilodeau and Srivastava (1992), and references therein). In one
of the first such studies, Muirhead and Verathaworn (1985) explain that the
ideal decision-theoretic approach that directly analyzes expected loss with
respect to the joint distribution of the eigenvalues of $B^{-1}A$
\textquotedblleft does not seem feasible due primarily to the complexity of
the distribution of the ordered latent roots...\textquotedblright\ Instead,
they focus on deriving an optimal estimator from a particular class.

Our \textit{LAN} result makes possible an asymptotic implementation of the
ideal decision-theoretic approach. We overcome the complexity of the joint
distribution of the eigenvalues by using a tractable contour integral
representation of the log likelihood process, which was obtained in the
single-spike case by Dharmawansa and Johnstone (2014). In the multiple-spike
case, a similar representation involves multiple contour integrals (see
Passemier et al (2014)). An asymptotic analysis of such a multiple integral
requires a substantial additional effort, and we leave it for future research.

It is interesting to contrast the \textit{LAN} result in the
\textit{super-critical} regime with the asymptotic behavior of the log
likelihood ratio in the case of a \textit{sub-critical} spike. In a separate
research effort, we follow Onatski et al (2013), who analyze the log
likelihood ratio in the \textit{sub-critical} regime for the case of a single
Wishart matrix, to show that the experiment of observing the eigenvalues of
$B^{-1}A$ in the \textit{sub-critical} regime is not of the \textit{LAN} type.
Furthermore, the log-likelihood process turns out to depend only on a smooth
functional of the empirical distribution of all the eigenvalues of
$\Sigma^{-1}A,$ so that asymptotically efficient inference procedures may
ignore the information contained in $\lambda_{1}$ altogether. The results of
this \textit{sub-critical} analysis will be published elsewhere.

The rest of the paper is structured as follows. In the next section, we
describe our setting. In Section 3, we derive the almost sure limits of a few
of the largest eigenvalues of the F-ratio. In Section 4, we establish the
asymptotic normality of the eigenvalue fluctuations in the
\textit{super-critical} regime. In Section 5, we derive an asymptotic
approximation to the joint distribution of the eigenvalues of $B^{-1}A$ for
the special case of a single \textit{super-critical} spike. In Section 6, we
show that the likelihood ratio in the local parameter space is asymptotically
equivalent to a centered and scaled largest eigenvalue, and establish the
\textit{LAN} property. Section 7 concludes.

\section{Setup}

Suppose that%
\[
A\sim W_{p}\left(  n_{1}+k,\Sigma_{1},\Omega_{1}\right)  \text{ \qquad and
}\qquad B\sim W_{p}\left(  n_{2},\Sigma_{2}\right)
\]
are independent non-central and central Wishart matrices respectively. For the
non-centrality parameter $\Omega_{1}$, we use a symmetric\ version of the
definition in Muirhead (1982, p. 442). That is, if $Z$ is an $n\times p$
matrix distributed as $N\left(  M,I_{n}\otimes\Sigma\right)  ,$ then
$Z^{\prime}Z\sim W_{p}\left(  n,\Sigma,\Omega\right)  $ with the
non-centrality parameter $\Omega=\Sigma^{-1/2}M^{\prime}M\Sigma^{-1/2}$. We
will consider two different settings for the parameters $\Sigma_{1},\Sigma
_{2},$ and $\Omega_{1}$.

\begin{description}
\item[Setting 1] \textbf{(Spiked covariance)} $\Sigma_{2}=\Sigma,$ $\Sigma
_{1}=\Sigma^{1/2}\left(  I_{p}+VhV^{\prime}\right)  \Sigma^{1/2},$ and
$\Omega_{1}=0.$ Here $\Sigma^{1/2}$ is the symmetric square root of a positive
definite matrix $\Sigma;$ $V$ in a $p\times k$ matrix of nuisance parameters
with orthonormal columns, and $h=\mathrm{diag}\left\{  h_{1},...,h_{k}%
\right\}  $ is the diagonal matrix of the \textquotedblleft covariance
spikes\textquotedblright\ with $h_{1}>...>h_{k}$.

\item[Setting 2] \textbf{(Spiked non-centrality)} $\Sigma_{2}=\Sigma,$
$\Sigma_{1}=\Sigma,$ and $\Omega_{1}=\left(  n_{1}+k\right)  VhV^{\prime},$
where $\Sigma,$ $V,$ and $h$ are as defined above, but $h_{j}$ with
$j=1,...,k$ are interpreted as \textquotedblleft non-centrality
spikes.\textquotedblright
\end{description}

We are interested in the behavior of the eigenvalues of
\[
\mathbf{F}\equiv\left(  B/n_{2}\right)  ^{-1}A/n_{A},
\]
where
\[
n_{A}=n_{1}+k,
\]
as $n_{1},n_{2},$ and $p$ grow so that $p/n_{1}\rightarrow c_{1}$ and
$p/n_{2}\rightarrow c_{2}$ with $0<c_{i}<1,$ while $k$, the number of spikes,
remains fixed. In what follows, we will assume that $\Sigma=I_{p}$. This
assumption is without loss of generality because the\ eigenvalues of
$\mathbf{F}$ do not change under the transformation $A\mapsto\Sigma
^{-1/2}A\Sigma^{-1/2},$ $B\mapsto\Sigma^{-1/2}B\Sigma^{-1/2}$.

It is convenient to think of $A/n_{A}$ as a sample covariance matrix
$XX^{\prime}/n_{A}$ of the sample $X$ having the factor structure%
\begin{equation}
X=V\mathcal{F}^{\prime}+\varepsilon\label{factor structure}%
\end{equation}
with $V,\mathcal{F},$ and $\varepsilon$ playing the roles of the factor
loadings, factors, and idiosyncratic terms, respectively. Matrices
$\mathcal{F}$ and $\varepsilon$ are mutually independent, and independent from
$B$. The distribution of $\varepsilon$ is $N\left(  0,I_{p}\otimes I_{n_{A}%
}\right)  ,$ and the distribution of $\mathcal{F}$ depends on the setting. For
Setting 1, $\mathcal{F}\sim N\left(  0,I_{p}\otimes h\right)  ,$ whereas for
Setting 2, $\mathcal{F}$ is a deterministic matrix such that $\mathcal{F}%
^{\prime}\mathcal{F}/n_{A}=h$. With this interpretation, Settings 1 and 2
describe, respectively, distributions which are unconditional and conditional
on the factors. In both cases the spike parameters $h_{j},$ $j=1,...,k,$
measure the factors' variability.

We would like to introduce a convenient representation of the eigenvalues of
$\mathbf{F}$, that we will denote as $\lambda_{p1}\geq...\geq\lambda_{pp}$.
First, note that $\lambda_{pj},$ $j=1,...,p$, are invariant with respect to
the simultaneous transformations
\begin{equation}
A\mapsto UAU^{\prime}\equiv n_{A}\tilde{H}\text{ \qquad and }\qquad B\mapsto
UBU^{\prime}\equiv n_{2}E, \label{random rotation}%
\end{equation}
where $U$ is a random matrix uniformly distributed over the orthogonal group
$\mathcal{O}\left(  p\right)  $. Under the assumption that $\Sigma=I_{p},$
matrix $n_{2}E$ is distributed as $W_{p}\left(  n_{2},I_{p}\right)  $ and is
independent from $\tilde{H}$. Matrix $\tilde{H}$ has the form $\tilde{X}%
\tilde{X}^{\prime}/n_{A},$ where%
\[
\tilde{X}=\tilde{V}\mathcal{F}^{\prime}+\tilde{\varepsilon}%
\]
with $\tilde{\varepsilon}\sim N\left(  0,I_{p}\otimes I_{n_{A}}\right)  $
independent from $\tilde{V}$, and $\tilde{V}$ being a random $p\times k$
matrix uniformly distributed on the Stiefel manifold of orthogonal $k$-frames
in $\mathbb{R}^{p}.$ We can think of $\tilde{V}$ as having the form
\[
\tilde{V}=v\left(  v^{\prime}v\right)  ^{-1/2}\equiv vW_{v}^{-1/2},
\]
where $v\sim N\left(  0,I_{p}\otimes I_{k}\right)  $ and $W_{v}\equiv
v^{\prime}v$ is Wishart $W_{k}\left(  p,I_{k}\right)  .$

Further, let $O_{\mathcal{F}}\in\mathcal{O}\left(  n_{A}\right)  $ be such
that the submatrix of its first $k$ columns equals $\mathcal{F}\left(
\mathcal{F}^{\prime}\mathcal{F}\right)  ^{-1/2}$, and let $\hat{X}=\tilde
{X}O_{\mathcal{F}}$. Clearly,
\begin{equation}
\tilde{H}=\tilde{X}\tilde{X}^{\prime}/n_{A}=\hat{X}\hat{X}^{\prime}/n_{A},
\label{H tilda}%
\end{equation}
and matrix $\hat{X}$ has the form%
\[
\hat{X}=vW_{v}^{-1/2}h^{1/2}W_{\mathcal{F}}^{1/2}+\hat{\varepsilon},
\]
where $v,W_{\mathcal{F}}$ and $\hat{\varepsilon}$ are mutually independent and
independent from $E;$ $\hat{\varepsilon}\sim N\left(  0,I_{p}\otimes I_{n_{A}%
}\right)  $; and the distribution of $W_{\mathcal{F}}$ depends on the setting.
For Setting~1, $W_{\mathcal{F}}\sim W_{k}\left(  n_{A},I_{k}\right)  ,$
whereas for Setting~2, $W_{\mathcal{F}}=n_{A}I_{k}$.

Finally, let us denote the submatrix of the first $k$ columns of
$\hat{\varepsilon}$ as $u.$ Then%
\begin{equation}
\hat{X}\hat{X}^{\prime}=\xi\xi^{\prime}+n_{1}H, \label{XhatXhat}%
\end{equation}
where $n_{1}H\sim W_{p}\left(  n_{1},I_{p}\right)  ,$ $H$ and $\xi\xi^{\prime
}$ are mutually independent, and independent from $E,$ and
\begin{equation}
\xi=vW_{v}^{-1/2}h^{1/2}W_{\mathcal{F}}^{1/2}+u. \label{perturbation vector}%
\end{equation}

Using (\ref{random rotation}), (\ref{H tilda}), and (\ref{XhatXhat}), we
obtain the convenient representation for the eigenvalues, announced above. Let
$\hat{x}_{p1}\geq...\geq\hat{x}_{pp}$ be the roots of the equation%
\begin{equation}
\det\left(  \xi\xi^{\prime}/n_{1}+H-xE\right)  =0. \label{canonical equation}%
\end{equation}
Then
\begin{equation}
\lambda_{pj}=n_{1}\hat{x}_{pj}/\left(  n_{1}+k\right)  .
\label{eigenvalues vs solutions}%
\end{equation}
This representation is convenient because the roots of
(\ref{canonical equation}) can be viewed and analyzed as perturbations of the
roots of equation $\det\left(  H-xE\right)  =0$ caused by adding the low-rank
matrix $\xi\xi^{\prime}/n_{1}$ to $H$.

If $x\in\mathbb{R}$ is such that $H-xE$ is invertible, then
\[
\left(  \xi\xi^{\prime}/n_{1}+H-xE\right)  ^{-1}=S-S\xi\left(  I_{k}%
+\xi^{\prime}S\xi/n_{1}\right)  ^{-1}\xi^{\prime}S/n_{1},
\]
where $S\equiv\left(  H-xE\right)  ^{-1}$. Therefore, if $x$ is a root of the
equation%
\begin{equation}
\det\left(  I_{k}+\xi^{\prime}\left(  H-xE\right)  ^{-1}\xi/n_{1}\right)  =0,
\label{equation for k largest}%
\end{equation}
then it also solves (\ref{canonical equation}), and hence, the asymptotic
behavior of the roots of (\ref{canonical equation}) can be inferred from that
of the random matrix-valued function%
\begin{equation}
M\left(  x\right)  =\xi^{\prime}\left(  H-xE\right)  ^{-1}\xi/n_{1}.
\label{quadratic form}%
\end{equation}
This is the main idea of the analysis in the next section of the paper.

\section{Almost sure limits of the largest eigenvalues}

Let $\mathbf{n}\equiv\left(  n_{1},n_{2}\right)  $ and $\mathbf{c}%
\equiv\left(  c_{1},c_{2}\right)  $. We will denote the asymptotic regime
where $n_{1},n_{2},$ and $p$ grow so that $p/n_{1}\rightarrow c_{1}$ and
$p/n_{2}\rightarrow c_{2}$ with $c_{j}\in\left(  0,1\right)  $ as
$p,\mathbf{n}\rightarrow_{\mathbf{c}}\infty$. As follows from Wachter's (1980)
work, as $p,\mathbf{n}\rightarrow_{\mathbf{c}}\infty$, the empirical
distribution of the eigenvalues of $E^{-1}H$ converges in probability to the
distribution with density%
\begin{equation}
\frac{1-c_{2}}{2\pi}\frac{\sqrt{\left(  b_{+}-\lambda\right)  \left(
\lambda-b_{-}\right)  }}{\lambda\left(  c_{1}+c_{2}\lambda\right)  }%
\mathbf{1}\left\{  b_{-}\leq\lambda\leq b_{+}\right\}  .
\label{Wachter density}%
\end{equation}
The upper and the lower boundaries of the support of this density are%
\[
b_{\pm}=\left(  \frac{1\pm r}{1-c_{2}}\right)  ^{2},\text{ where }%
r=\sqrt{c_{1}+c_{2}-c_{1}c_{2}}\text{.}%
\]
The results of Silverstein and Bai (1995) and Silverstein (1995) show that the
empirical distribution converges not only in probability, but also almost
surely (a.s.). Furthermore, as follows from Theorem 1.1 of Bai and Silverstein
(1998), the largest eigenvalue of $E^{-1}H$ a.s. converges to $b_{+}$.

The latter convergence, together with (\ref{eigenvalues vs solutions}) and
Weyl's inequalities for the eigenvalues of a sum of two Hermitian matrices
(see Theorem 4.3.7 in Horn and Johnson (1985)), imply that the $k+1$-th
largest eigenvalue of $\mathbf{F,}$ $\lambda_{p,k+1}$, a.s. converges to
$b_{+}$. Those of the $k$ largest eigenvalues that remain separated from
$b_{+}$ as $p,\mathbf{n}\rightarrow_{\mathbf{c}}\infty$, must correspond to
solutions of (\ref{equation for k largest}). Below, we study these solutions
in detail.

\begin{lemma}
\label{Lemma1}For any $x>b_{+},$ as $p,\mathbf{n}\rightarrow_{\mathbf{c}%
}\infty$,%
\begin{align}
&  \frac{1}{p}\operatorname*{tr}\left[  \left(  H-xE\right)  ^{-1}\right]
\overset{a.s.}{\rightarrow}m_{x}(0)\text{ and}\label{as lim}\\
&  \frac{1}{p}\operatorname*{tr}\left[  \frac{\mathrm{d}}{\mathrm{d}x}\left(
H-xE\right)  ^{-1}\right]  \overset{a.s.}{\rightarrow}\frac{\mathrm{d}%
}{\mathrm{d}x}m_{x}(0), \label{as lim dx}%
\end{align}
where $m_{x}(0)=\lim_{z\rightarrow0}m_{x}(z),$ and $m_{x}(z)\in\mathbb{C}^{+}$
is analytic in $z\in\mathbb{C}^{+},$ and satisfies equation%
\begin{equation}
z-\frac{1}{1+c_{1}m_{x}\left(  z\right)  }=-\frac{1}{m_{x}\left(  z\right)
}-\frac{x}{1-c_{2}xm_{x}\left(  z\right)  }. \label{eq3}%
\end{equation}

\end{lemma}

\textbf{Proof:} Let $x\in\mathbb{R}$ be such that $x>b_{+},$ and let
$F_{x}(\lambda)$ be the empirical distribution function of the eigenvalues of
$H-xE$. For any $z\in\mathbb{C}^{+},$ let
\[
\hat{m}_{x}\left(  z\right)  =\int\left(  \lambda-z\right)  ^{-1}%
\mathrm{d}F_{x}(\lambda)
\]
be the Stieltjes transform of $F_{x}(\lambda)$. Note that matrix $H-xE$ can be
represented in the form $YTY^{\prime}/p,$ where $Y\sim N\left(  0,I_{p}\otimes
I_{n_{1}+n_{2}}\right)  $ and $T$ is a diagonal matrix with the first $n_{1}$
and the last $n_{2}$ diagonal elements equal to $p/n_{1}$ and $-xp/n_{2},$
respectively. Therefore, by Theorem 1.1 of Silverstein and Bai (1995), for any
$z\in\mathbb{C}^{+},$ $\hat{m}_{x}\left(  z\right)  $ a.s. converges to
$m_{x}\left(  z\right)  \in\mathbb{C}^{+},$ which is an analytic function in
the domain $z\in\mathbb{C}^{+}$ that solves the functional equation (\ref{eq3}).

By Theorem 1.1 of Bai and Silverstein (1998), the largest eigenvalue of
$E^{-1}H$ a.s. converges to $b_{+}$. Therefore, for any $x>b_{+},$ the largest
eigenvalue of $H-xE$ is a.s. asymptotically bounded away from the positive
semi-axis. Hence, $\hat{m}_{x}\left(  z\right)  $ is analytic and bounded in a
small disc $D$ around $z=0$ for all sufficiently large $p$ and $\mathbf{n},$
a.s. By Vitali's theorem (see Titchmarsh (1960), p.168), $\hat{m}_{x}\left(
z\right)  $ is a.s. converging to an analytic function in $D$. Since, in
$D\cap\mathbb{C}^{+}$, the limiting function is $m_{x}\left(  z\right)  ,$ we
have
\[
\frac{1}{p}\operatorname*{tr}\left[  \left(  H-xE\right)  ^{-1}\right]
=\hat{m}_{x}\left(  0\right)  \overset{a.s.}{\rightarrow}m_{x}\left(
0\right)  ,
\]
where $m_{x}(0)=\lim_{z\rightarrow0}m_{x}(z)$. Further, $\frac{1}%
{p}\operatorname*{tr}\left[  \left(  H-\zeta E\right)  ^{-1}\right]  $ is an
analytic bounded function of $\zeta$ in a small disk $D_{x}$ around $x,$ for
all sufficiently large $p$ and $\mathbf{n},$ a.s. Therefore, by Vitali's
theorem its a.s. limit $f(\zeta)$ is analytic in $D_{x},$ and
\[
\frac{1}{p}\operatorname*{tr}\left[  \frac{\mathrm{d}}{\mathrm{d}\zeta}\left(
H-\zeta E\right)  ^{-1}\right]  \overset{a.s.}{\rightarrow}\frac{\mathrm{d}%
}{\mathrm{d}\zeta}f(\zeta)
\]
in $D_{x}.$ On the other hand, we know that $f(\zeta)=m_{\operatorname{Re}%
\zeta}(0)$ for $\zeta$ from $D_{x}$. Therefore, we have (\ref{as lim dx}%
).$\square$

\begin{lemma}
\label{Lemma2}For any $x>b_{+}$, as $p,\mathbf{n}\rightarrow_{\mathbf{c}%
}\infty$,%
\begin{align*}
&  \left\Vert M\left(  x\right)  -\left(  h+c_{1}I_{k}\right)  \frac{1}%
{p}\operatorname*{tr}\left[  \left(  H-xE\right)  ^{-1}\right]  \right\Vert
\overset{a.s.}{\rightarrow}0\text{ and}\\
&  \left\Vert \frac{\mathrm{d}}{\mathrm{d}x}M\left(  x\right)  -\left(
h+c_{1}I_{k}\right)  \frac{1}{p}\operatorname*{tr}\left[  \frac{\mathrm{d}%
}{\mathrm{d}x}\left(  H-xE\right)  ^{-1}\right]  \right\Vert
\overset{a.s.}{\rightarrow}0,
\end{align*}
where $\left\Vert \cdot\right\Vert $ denotes the spectral norm.
\end{lemma}

\textbf{Proof:} This convergences follow from (\ref{perturbation vector}),
(\ref{quadratic form}), and Lemma \ref{Bai-Silverstein}\ stated below.$\square
$

\begin{lemma}
\label{Bai-Silverstein} Let $C$ be a random $p\times p$ matrix, independent
from $u$ and $v$, which are as defined in Section 2, and such that
$p\left\Vert C\right\Vert $ is bounded for all sufficiently large $p,$ a.s.
Then, as $p\rightarrow\infty$,%
\[
\left\Vert v^{\prime}Cv-\left(  \operatorname*{tr}C\right)  I_{k}\right\Vert
\overset{a.s.}{\rightarrow}0\text{ and }\left\Vert v^{\prime}Cu\right\Vert
\overset{a.s.}{\rightarrow}0.
\]

\end{lemma}

\textbf{Proof:} This lemma follows from the Borel-Cantelli lemma, and the
upper bounds on the fourth moments of the entries $v^{\prime}Cv-\left(
\operatorname*{tr}C\right)  I_{k}$ and $v^{\prime}Cu$ established by Lemma 2.7
of Bai and Silverstein (1998).$\square$

\begin{lemma}
\label{monotonicity} (i) For any $\varepsilon>0,$ the $k$ eigenvalues of
$M\left(  x\right)  $ are strictly increasing functions of $x\in\left(
b_{+}+\varepsilon,\infty\right)  $ for sufficiently large $p$ and $\mathbf{n}%
$, a.s.; (ii) $m_{x}(0)$ is a strictly increasing, continuous function of
$x\in\left(  b_{+},\infty\right)  $; (iii) $\lim_{x\rightarrow\infty}%
m_{x}(0)=0$, and $\lim_{x\downarrow b_{+}}m_{x}(0)\left(  h_{i}+c_{1}\right)
<-1$ if and only if $h_{i}>\bar{h},$ where%
\[
\bar{h}=\frac{c_{2}+r}{1-c_{2}}.
\]

\end{lemma}

\textbf{Proof:} Let $\mu_{1}\in\left(  0,\infty\right)  $ be the largest
eigenvalue of $E^{-1}H.$ For any $x_{1}>x_{2}>\mu_{1},$ matrix $\left(
x_{1}E-H\right)  ^{-1}-\left(  x_{2}E-H\right)  ^{-1}$ is negative definite,
a.s. Part (i) follows from this, from the definition (\ref{quadratic form}),
and from the fact that $\mu_{1}\overset{a.s.}{\rightarrow}b_{+}.$ Part (i)
together with Lemmas \ref{Lemma1} and \ref{Lemma2} imply that $m_{x}(0)$ is
increasing on $\left(  b_{+},\infty\right)  .$ It is strictly increasing
because, otherwise, (\ref{eq3}) would not be satisfied for some $z\in
\mathbb{C}^{+}$ that are sufficiently close to zero. The continuity follows
from the analyticity of $m_{x}(0)$ established in the proof of Lemma
\ref{Lemma1}. Finally, $\lim_{x\rightarrow\infty}m_{x}(0)=0$ is implied by
(ii) and (\ref{as lim}). Equation (\ref{eq3}) implies that
\[
\lim_{x\downarrow b_{+}}m_{x}(0)=\frac{c_{2}-1}{\left(  r+1\right)  r},
\]
which, in its turn, implies the second statement of (iii).$\square$

Let $\hat{x}_{p1}\geq...\geq\hat{x}_{pk}$ be the solutions of equation
(\ref{equation for k largest}). By Lemmas \ref{Lemma1}, \ref{Lemma2}, and
\ref{monotonicity}, if
\begin{equation}
h_{1}>...>h_{m}>\bar{h}>h_{m+1}>...>h_{k}, \label{general h}%
\end{equation}
then $\hat{x}_{pi}\overset{a.s.}{\rightarrow}x_{i},$ where $x_{i},$
$i=1,...,m,$ are such that
\begin{equation}
1+\left(  h_{i}+c_{1}\right)  m_{x_{i}}(0)=0 \label{root}%
\end{equation}
and $m_{x_{i}}(0)$ satisfies (\ref{eq3}) with $x$ replaced by $x_{i}.$ In
particular,%
\begin{equation}
\frac{1}{1+c_{1}m_{x_{i}}(0)}-\frac{1}{m_{x_{i}}(0)}-\frac{x_{i}}{1-c_{2}%
x_{i}m_{x_{i}}(0)}=0. \label{implicit}%
\end{equation}

Combining (\ref{root}) and (\ref{implicit}), we obtain%
\[
\frac{1}{h_{i}}+1-\frac{x_{i}}{h_{i}+c_{1}+c_{2}x_{i}}=0,
\]
which implies that%
\begin{equation}
x_{i}=\frac{\left(  h_{i}+c_{1}\right)  \left(  h_{i}+1\right)  }{h_{i}%
-c_{2}\left(  h_{i}+1\right)  }. \label{x0}%
\end{equation}
By (\ref{eigenvalues vs solutions}), $n_{1}\hat{x}_{pi}/\left(  n_{1}%
+k\right)  ,$ $i=1,...,m$, must be the $m$ largest eigenvalues of
$\mathbf{F},$ and thus, $x_{i},$ $i=1,...,m,$ describe their a.s. limits.
Since there are only $m$ roots of (\ref{equation for k largest}) that are
asymptotically separated from $b_{+}$ and are located above $b_{+},$ the other
$k-r$ of the largest eigenvalues of $\mathbf{F}$ must a.s. converge to $b_{+}%
$. To summarize, the following proposition holds.

\begin{proposition}
\label{Proposition1}Suppose that $h_{1}>...>h_{m}>\bar{h}>h_{m+1}>...>h_{k}$.
Then, for $i\leq m,$ the $i$-th largest eigenvalue of $\mathbf{F}$ a.s.
converges to $x_{i}$ defined in (\ref{x0}). For $m<i\leq k,$ the $i$-th
largest eigenvalue a.s. converges to $b_{+}.$
\end{proposition}

As follows from Proposition \ref{Proposition1}, $\bar{h}=\left(
c_{2}+r\right)  /\left(  1-c_{2}\right)  $ is the phase transition threshold
for the eigenvalues of the spiked F-ratio $\mathbf{F}$. The value of this
threshold diverges to infinity when $c_{2}\rightarrow1$. Note that, when
$c_{2}$ is close to one, the smallest eigenvalue of $B/n_{2}$ is close to
zero, which makes $\left(  B/n_{2}\right)  ^{-1}$ a particularly bad estimator
of the inverse of the population covariance, $\Sigma^{-1}$. When
$c_{2}\rightarrow0,$ the phase transition converges to $\sqrt{c_{1}},$ which
is the phase transition threshold for the eigenvalues of a single spiked
Wishart matrix. In such a case, $x_{i}$ converges to $\left(  h_{i}%
+c_{1}\right)  \left(  h_{i}+1\right)  /h_{i},$ which is the a.s. limit of the
$i$-th largest eigenvalue of the spiked Wishart when the $i$-th spike $h_{i}$
is above $\sqrt{c_{1}}$.

\section{Asymptotic normality}

In what follows, we will assume that (\ref{general h}) holds, so that only $m$
eigenvalues of $\mathbf{F}$ separate from the bulk asymptotically. We would
like to study their fluctuations around the corresponding a.s. limits.
Proposition \ref{Proposition1} shows that the limits $x_{i}$ depend on $c_{1}$
and $c_{2}$. Because of this dependence, the rate of the convergence has to
depend on the rates of the convergences $p/n_{1}\rightarrow c_{1}$ and
$p/n_{2}\rightarrow c_{2}$. However, as will be shown below, the latter rates
do not affect the fluctuations of $\lambda_{pi}$ around%
\[
x_{pi}=\frac{\left(  h_{i}+c_{p1}\right)  \left(  h_{i}+1\right)  }%
{h_{i}-c_{p2}\left(  h_{i}+1\right)  },
\]
which are obtained from $x_{i}$ by replacing $c_{1}$ and $c_{2}$ by
$c_{p1}=p/n_{1}$ and $c_{p2}=p/n_{2}$.

Similar to $x_{i}$, which are linked to the Stieltjes transform of the
limiting spectral distribution of $xE-H$ via (\ref{root}), $x_{pi}$ also can
be linked to the limiting Stieltjes transform, albeit under a slightly
different asymptotic regime. Precisely, let $m_{px}\left(  z\right)  $ be the
Stieltjes transform of the limiting spectral distribution of $xE-H$ as
$n_{1},$ $n_{2},$ and $p$ grow so that $p/n_{1}$ and $p/n_{2}$ \textit{remain
fixed}. Then, similarly to (\ref{root}), we have%
\begin{equation}
1+\left(  h_{i}+c_{p1}\right)  m_{px_{pi}}(0)=0. \label{rootp}%
\end{equation}
This link will be useful in our analysis below, where we maintain the
assumption that $p/n_{1}$ and $p/n_{2}$ are not necessarily fixed, but
converge to $c_{1}$ and $c_{2},$ respectively.

Recall that, by (\ref{eigenvalues vs solutions}), $\lambda_{pi}=n_{1}\hat
{x}_{pi}/\left(  n_{1}+k\right)  ,$ where $\hat{x}_{pi},$ $i=1,...,m,$ satisfy
(\textbf{\ref{equation for k largest}}). Clearly, the asymptotic distributions
of $\sqrt{p}\left(  \lambda_{pi}-x_{pi}\right)  $ and $\sqrt{p}\left(  \hat
{x}_{pi}-x_{pi}\right)  ,$ $i=1,...,m,$ coincide. Therefore, below we will
study the asymptotic behavior of $\sqrt{p}\left(  \hat{x}_{pi}-x_{pi}\right)
,$ $i=1,...,m.$ By the standard Taylor expansion argument,%
\begin{equation}
\sqrt{p}\left(  \hat{x}_{pi}-x_{pi}\right)  =-\frac{\sqrt{p}\det
\mathcal{M}\left(  x_{pi}\right)  }{\frac{\mathrm{d}}{\mathrm{d}x}%
\det\mathcal{M}\left(  x_{pi}\right)  +\frac{1}{2}\left(  \hat{x}_{pi}%
-x_{pi}\right)  \frac{\mathrm{d}^{2}}{\mathrm{d}x^{2}}\det\mathcal{M}\left(
\tilde{x}_{pi}\right)  }, \label{Taylor argument}%
\end{equation}
$i=1,...,m$, where $\mathcal{M}\left(  x\right)  =I_{k}+M\left(  x\right)  ,$
and $\tilde{x}_{pi}\in\left[  x_{pi},\hat{x}_{pi}\right]  .$ We have%
\[
\frac{\mathrm{d}}{\mathrm{d}x}\det\mathcal{M}\left(  x_{pi}\right)
=\det\mathcal{M}\left(  x_{pi}\right)  \operatorname*{tr}S\left(
x_{pi}\right)  ,
\]
and%
\[
\frac{\mathrm{d}^{2}}{\mathrm{d}x^{2}}\det\mathcal{M}\left(  x_{pi}\right)
=\det\mathcal{M}\left(  x_{pi}\right)  \left\{  \operatorname*{tr}R\left(
x_{pi}\right)  +\left(  \operatorname*{tr}S\left(  x_{pi}\right)  \right)
^{2}-\operatorname*{tr}\left[  S^{2}\left(  x_{pi}\right)  \right]  \right\}
,
\]
where%
\[
S(x)=\mathcal{M}\left(  x\right)  ^{-1}\frac{\mathrm{d}}{\mathrm{d}%
x}M(x),\text{ and }R(x)=\mathcal{M}\left(  x\right)  ^{-1}\frac{\mathrm{d}%
^{2}}{\mathrm{d}x^{2}}M(x).
\]
Since the event%
\[
\det\mathcal{M}\left(  x_{pi}\right)  =0\text{ or }1+M_{ii}(x_{pi})=0\text{
for some }i=1,...,m
\]
happens with probability zero, we can simultaneously multiply the numerator
and denominator of (\ref{Taylor argument}) by $\left(  1+M_{ii}(x_{pi}%
)\right)  /\det\mathcal{M}\left(  x_{pi}\right)  $ to obtain%
\begin{equation}
\sqrt{p}\left(  \hat{x}_{pi}-x_{pi}\right)  =-\frac{\sqrt{p}\left(
1+M_{ii}(x_{pi})\right)  }{s(x_{pi})+\frac{1}{2}\left(  \hat{x}_{pi}%
-x_{pi}\right)  \delta(x_{pi})}, \label{Taylor1}%
\end{equation}
where%
\[
s(x_{pi})=\left(  1+M_{ii}(x_{pi})\right)  \operatorname*{tr}S\left(
x_{pi}\right)  ,
\]
and%
\[
\delta(x_{pi})=\left(  1+M_{ii}(x_{pi})\right)  \left\{  \operatorname*{tr}%
R\left(  x_{pi}\right)  +\left(  \operatorname*{tr}S\left(  x_{pi}\right)
\right)  ^{2}-\operatorname*{tr}\left[  S^{2}\left(  x_{pi}\right)  \right]
\right\}  .
\]

\begin{lemma}
\label{denominator}For any $i=1,...,m,$ we have: (i) $s(x_{pi}%
)\overset{\mathrm{P}}{\rightarrow}\left(  h_{i}+c_{1}\right)  \frac
{\mathrm{d}}{\mathrm{d}x}m_{x_{i}}(0)$; (ii)~$\delta(x_{pi})=O(1)$ a.s.
\end{lemma}

\textbf{Proof}: By Lemmas \ref{Lemma1} and \ref{Lemma2},
\begin{equation}
\frac{\mathrm{d}}{\mathrm{d}x}M(x_{pi})\overset{a.s.}{\rightarrow}\left(
h+c_{1}I_{k}\right)  \frac{\mathrm{d}}{\mathrm{d}x}m_{x_{i}}(0).
\label{simple convergence}%
\end{equation}
Further,
\begin{equation}
\left(  1+M_{ii}(x_{pi})\right)  \left(  I_{k}+M(x_{pi})\right)
^{-1}\overset{a.s.}{\rightarrow}\operatorname*{diag}\left\{
0,...,0,1,0,...,0\right\}  \label{ii convergence}%
\end{equation}
with 1 at the $i$-th place on the diagonal. The latter convergence follows
from the fact that $I_{k}+M(x_{pi})$ can be viewed as a small perturbation of
a diagonal matrix
\[
I_{k}+\left(  h+c_{1}I_{k}\right)  m_{x_{i}}(0),
\]
which has non-zero diagonal elements, except at the $i$-th position. The
eigenvalue perturbation formulae (see, for example, (2.33) on p.79 of Kato
(1980)) will then lead to (\ref{ii convergence}). Combining
(\ref{simple convergence}) and (\ref{ii convergence}), and using the
definition of $s(x_{pi}),$ we obtain (i).

To establish (ii), we note that $\left(  1+M_{ii}(x_{pi})\right)
\operatorname*{tr}R\left(  x_{pi}\right)  =O_{\mathrm{P}}(1)$ by an argument
similar to that used to establish (i). Further, $\left(  \operatorname*{tr}%
S\left(  x_{pi}\right)  \right)  ^{2}-\operatorname*{tr}\left[  S^{2}\left(
x_{pi}\right)  \right]  $ is a linear function of the only eigenvalue of
$S\left(  x_{pi}\right)  $ that diverges to infinity. By the eigenvalue
perturbation formulae, such an eigenvalue equals $\left(  1+M_{ii}%
(x_{pi})\right)  ^{-1}O(1)$ a.s. Therefore,%
\[
\left(  1+M_{ii}(x_{pi})\right)  \left(  \left(  \operatorname*{tr}S\left(
x_{pi}\right)  \right)  ^{2}-\operatorname*{tr}\left[  S^{2}\left(
x_{pi}\right)  \right]  \right)  =O(1),
\]
which concludes the proof of (ii).$\square$

Equation (\ref{Taylor1}), Lemma \ref{denominator}, and the Slutsky theorem
imply that, for the purpose of establishing convergence in distribution of
$\sqrt{p}\left(  \hat{x}_{pi}-x_{pi}\right)  $, $i=1,...,m,$ we may focus on
the numerator of (\ref{Taylor1})%
\[
Z_{ii}(x_{pi})\equiv\sqrt{p}\left(  1+M_{ii}(x_{pi})\right)  =\sqrt{p}\left(
M_{ii}(x_{pi})-\left(  h_{i}+c_{p1}\right)  m_{px_{pi}}(0)\right)  ,
\]
where the last equality follows from (\ref{rootp}).

The random variable $Z_{ii}$ is the entry of the matrix%
\[
Z(x_{pi})=\sqrt{p}\left(  M(x_{pi})-\left(  h+c_{p1}I_{k}\right)  m_{px_{pi}%
}(0)\right)
\]
that belongs to the $i$-th row and the $i$-th column. Let us now introduce new
notations. Let%
\begin{align*}
D  &  =\left(  W_{\mathcal{F}}/n_{1}\right)  ^{1/2}h^{1/2}\left(
W_{v}/p\right)  ^{-1/2},\\
G  &  =\left(  H-x_{pi}E\right)  ^{-1}/p,\\
\Delta_{\mathcal{F}}  &  =\sqrt{n_{1}}\left(  \left(  W_{\mathcal{F}}%
/n_{1}\right)  ^{1/2}-I_{k}\right)  ,\text{ and}\\
\Delta_{v}  &  =\sqrt{p}\left(  W_{v}/p-I_{k}\right)  .
\end{align*}
Then, using equations (\ref{quadratic form}) and (\ref{perturbation vector}),
we obtain the following decomposition.%
\[
Z(x_{pi})=%
{\displaystyle\sum\nolimits_{v=1}^{6}}
Z^{(v)},
\]
where%
\[
Z^{(1)}=D\sqrt{p}\left(  v^{\prime}Gv-I_{k}\operatorname*{tr}G\right)
D^{\prime},
\]%
\[
Z^{(2)}=\left(  \operatorname*{tr}G\right)  D\left(  W_{v}/p\right)
^{-1/2}h^{1/2}\sqrt{c_{p1}}\Delta_{\mathcal{F}},
\]%
\[
Z^{(3)}=\operatorname*{tr}G\sqrt{c_{p1}}\Delta_{\mathcal{F}}h^{1/2}\left(
W_{v}/p\right)  ^{-1}h^{1/2},
\]%
\[
Z^{(4)}=-\left(  \operatorname*{tr}G\right)  h^{1/2}\Delta_{v}\left(
W_{v}/p\right)  ^{-1}h^{1/2},
\]%
\[
Z^{(5)}=\sqrt{c_{p1}}\sqrt{p}\left(  Dv^{\prime}Gu+u^{\prime}GvD^{\prime
}\right)  ,
\]%
\[
Z^{(6)}=c_{p1}\sqrt{p}\left(  u^{\prime}Gu-I_{k}\operatorname*{tr}G\right)  ,
\]
and%
\[
Z^{(7)}=\left(  h+c_{p1}I_{k}\right)  \sqrt{p}\left(  \operatorname*{tr}%
G-m_{px_{pi}}(0)\right)  .
\]
For the last term, $Z^{(7)},$ we prove the following lemma.

\begin{lemma}
\label{speed}$Z^{(7)}\overset{\mathrm{P}}{\rightarrow}0.$
\end{lemma}

\textbf{Proof:} The proof of this lemma will appear in a separate work. Had
$x_{pi}$ been negative, $H-x_{pi}E$ would have been having the form
$YTY^{\prime}$ with $Y\sim N\left(  0,I_{p}\otimes I_{n_{1}+n_{2}}\right)  $
and a positive definite diagonal $T$ with converging spectral distribution.
The Lemma would have been following then from the results of Bai and
Silverstein (2004). Our proof extends Bai and Silverstein's (2004) arguments
to the case of negative $x_{pi}$.$\square$

Further, the asymptotic behavior of the terms $Z^{(2)}$ and $Z^{(3)}$ differ
depending on the setting. Recall that for Setting 1, $W_{\mathcal{F}}\sim
W_{k}\left(  n_{A},I_{k}\right)  $. Then, since%
\[
\Delta_{\mathcal{F}}=\sqrt{n_{1}}\left(  W_{\mathcal{F}}/n_{1}-I_{k}\right)
/2+o_{\mathrm{P}}\left(  1\right)  ,
\]
a standard CLT together with Lemma \ref{Lemma1} imply that%
\begin{equation}
Z^{(2)}+Z^{(3)}\overset{d}{\rightarrow}N\left(  0,2c_{1}m_{x_{i}}^{2}%
(0)h^{2}\right)  . \label{Z2}%
\end{equation}
The latter limit is independent from the limits of $Z^{(j)},$ $j\neq2,3,$
because $W_{\mathcal{F}}$ is independent from $u$ and $v$.

In contrast, for Setting 2, we have $W_{\mathcal{F}}=n_{A}I_{k},$ and
$\Delta_{\mathcal{F}}=o(1).$ Therefore,%
\begin{equation}
Z^{(2)}+Z^{(3)}\overset{\mathrm{P}}{\rightarrow}0. \label{Z2 setting 2}%
\end{equation}

Let us now establish the convergence of $Z^{(j)},$ $j\leq6$ such that
$j\neq2,3$. Let $l_{i}$ and $L_{i}$ be such that $\left[  l_{i},L_{i}\right]
$ includes the support of the limiting spectral distribution, $G_{x_{i}}$, of
$H-x_{pi}E$. Moreover, let $\left[  l_{i},L_{i}\right]  $ be such that none of
the eigenvalues $\lambda_{p1}^{(i)},...,\lambda_{pp}^{(i)}$ of $H-x_{pi}E$
lies outside $\left[  l_{i},L_{i}\right]  $ for sufficiently large $p$, a.s.
Further, let $g_{q}$ with $q=1,...,Q,$ where $Q$ is an arbitrary positive
integer, be functions which are continuous on $\left[  l_{i},L_{i}\right]  $
and let $\zeta$ denote a $p\times m$ matrix with i.i.d. $N(0,1)$ entries$.$
Finally, let
\[
\Theta=\left\{  \left(  q,s,t\right)  :q=1,...,Q;1\leq s\leq t\leq m\right\}
.
\]
The following Lemma is a slight modification of Lemma 13 of the Supplementary
Appendix in Onatski (2012).

\begin{lemma}
\label{Lemma 13}The joint distribution of random variables%
\[
\left\{  \frac{1}{\sqrt{p}}%
{\displaystyle\sum\nolimits_{j=1}^{p}}
g_{q}\left(  \lambda_{pj}^{(i)}\right)  \left(  \zeta_{js}\zeta_{jt}%
-\delta_{st}\right)  ,\left(  q,s,t\right)  \in\Theta\right\}
\]
weakly converges to a multivariate normal. The covariance between components
$\left(  q,s,t\right)  $ and $\left(  q_{1},s_{1},t_{1}\right)  $ of the
limiting distribution is equal to $0$ when $\left(  s,t\right)  \neq\left(
s_{1},t_{1}\right)  ,$ and to $\left(  1+\delta_{st}\right)  \int
g_{q}(\lambda)g_{q_{1}}\left(  \lambda\right)  \mathrm{d}G_{x_{i}}(\lambda)$
when $\left(  s,t\right)  =\left(  s_{1},t_{1}\right)  $.
\end{lemma}

\textbf{Proof:} For readers' convenience, we provide a proof of this Lemma in
the Appendix.$\square$

Note that all entries of $Z^{(j)},$ $j\leq6$ such that $j\neq2,3,$ are linear
combinations of the terms having the form considered in Lemma \ref{Lemma 13},
with weights converging in probability to finite constants. Take, for example
$Z^{(1)}$. Its entries are linear combinations of the entries of
\[
\frac{1}{\sqrt{p}}v^{\prime}\left(  H-x_{pi}E\right)  ^{-1}v-I_{k}\frac
{1}{\sqrt{p}}\operatorname*{tr}\left(  H-x_{pi}E\right)  ^{-1},
\]
which, in turn, can be represented in the form $\frac{1}{\sqrt{p}}%
{\displaystyle\sum\nolimits_{j=1}^{p}}
\left(  \lambda_{pj}^{(i)}\right)  ^{-1}\left(  \zeta_{js}\zeta_{jt}%
-\delta_{st}\right)  .$ The matrix $\zeta$ is obtained by multiplying $\left[
u,v\right]  $ from the left by the eigenvector matrix of $H-x_{pi}E$.

Lemma \ref{Lemma 13} implies that vector $\left(  Z_{ii}^{(1)},Z_{ii}%
^{(4)},Z_{ii}^{(5)},Z_{ii}^{(6)}\right)  $ converges in distribution to a
four-dimensional normal vector with zero mean and the following covariance
matrix%
\[
\left(
\begin{array}
[c]{cccc}%
2h_{i}^{2}m_{x_{i}}^{\prime}\left(  0\right)  & -2h_{i}^{2}m_{x_{i}}%
^{2}\left(  0\right)  & 0 & 0\\
-2h_{i}^{2}m_{x_{i}}^{2}\left(  0\right)  & 2h_{i}^{2}m_{x_{i}}^{2}\left(
0\right)  & 0 & 0\\
0 & 0 & 4c_{1}h_{i}m_{x_{i}}^{\prime}\left(  0\right)  & 0\\
0 & 0 & 0 & 2c_{1}^{2}m_{x_{i}}^{\prime}\left(  0\right)
\end{array}
\right)  .
\]
Combining this result with Lemma \ref{speed}, and convergencies (\ref{Z2}),
and (\ref{Z2 setting 2}), we obtain, for Setting 1,
\begin{equation}
Z_{ii}(x_{pi})\overset{d}{\rightarrow}N\left(  0,2\left(  h_{i}+c_{1}\right)
^{2}m_{x_{i}}^{\prime}\left(  0\right)  -2h_{i}^{2}\left(  1-c_{1}\right)
m_{x_{i}}^{2}(0)\right)  , \label{convergence in distribution}%
\end{equation}
and, for Setting 2,%
\begin{equation}
Z_{ii}(x_{pi})\overset{d}{\rightarrow}N\left(  0,2\left(  h_{i}+c_{1}\right)
^{2}m_{x_{i}}^{\prime}\left(  0\right)  -2h_{i}^{2}m_{x_{i}}^{2}(0)\right)  .
\label{convergence in d set 2}%
\end{equation}
To establish the joint convergence of $Z_{ii}(x_{pi}),$ $i=1,...,m,$ we need
another lemma. For each $i=1,...,m,$ let $g_{q}^{(i)},$ with $q=1,...,Q,$ be
functions continuous on $\left[  l_{i},L_{i}\right]  .$

\begin{lemma}
\label{Lemma 13a}For any set of pairs $\left\{  \left(  s_{i},t_{i}\right)
:i=1,...,m\right\}  $ such that $\left(  s_{i_{1}},t_{i_{1}}\right)
\neq\left(  s_{i_{2}},t_{i_{2}}\right)  $ for any $i_{1}\neq i_{2},$ the joint
distribution of random variables%
\[
\left\{  \frac{1}{\sqrt{p}}%
{\displaystyle\sum\nolimits_{j=1}^{p}}
g_{q}^{(i)}\left(  \lambda_{pj}^{(i)}\right)  \left(  \zeta_{js_{i}}%
\zeta_{jt_{i}}-\delta_{s_{i}t_{i}}\right)  ,i=1,...,m\right\}
\]
weakly converges to a multivariate normal. The covariance between components
$i_{1}$ and $i_{2}$ of the limiting distribution is equal to $0$ when
$i_{1}\neq i_{2}.$
\end{lemma}

The proof of this lemma is very similar to that of Lemma \ref{Lemma 13}, and
we omit it to save space. Lemma \ref{Lemma 13a} implies that $Z_{ii}(x_{pi}),$
$i=1,...,m$ jointly converge to an $m$-dimensional normal vector with a
diagonal covariance matrix. This result, together with equation (\ref{Taylor1}%
), Lemma \ref{denominator}, and convergences
(\ref{convergence in distribution}, \ref{convergence in d set 2}) establish
the following Lemma.

\begin{lemma}
The joint asymptotic distribution of $\sqrt{p}\left(  \lambda_{pi}%
-x_{pi}\right)  ,i=1,...,m$ is normal, with diagonal covariance matrix. For
Setting 1, the $i$-th diagonal element of the covariance matrix equals
\begin{equation}
\frac{2\left(  h_{i}+c_{1}\right)  ^{2}m_{x_{i}}^{\prime}\left(  0\right)
-2h_{i}^{2}\left(  1-c_{1}\right)  m_{x_{i}}^{2}(0)}{\left(  h_{i}%
+c_{1}\right)  ^{2}\left(  \frac{\mathrm{d}}{\mathrm{d}x}m_{x_{i}}(0)\right)
^{2}}. \label{normality}%
\end{equation}
For Setting 2, it equals%
\begin{equation}
\frac{2\left(  h_{i}+c_{1}\right)  ^{2}m_{x_{i}}^{\prime}\left(  0\right)
-2h_{i}^{2}m_{x_{i}}^{2}(0)}{\left(  h_{i}+c_{1}\right)  ^{2}\left(
\frac{\mathrm{d}}{\mathrm{d}x}m_{x_{i}}(0)\right)  ^{2}}.
\label{normality Setting 2}%
\end{equation}

\end{lemma}

In the Appendix, we establish the following explicit expressions for
$m_{x_{i}}^{2}\left(  0\right)  ,$ $m_{x_{i}}^{\prime}\left(  0\right)  ,$ and
$\frac{\mathrm{d}}{\mathrm{d}x}m_{x_{i}}(0):$
\begin{equation}
m_{x_{i}}^{2}\left(  0\right)  =\left(  h_{i}+c_{1}\right)  ^{-2},
\label{ingredient1}%
\end{equation}%
\begin{equation}
m_{x_{i}}^{\prime}\left(  0\right)  =-\frac{h_{i}^{2}}{\left(  h_{i}%
+c_{1}\right)  ^{2}\left(  c_{1}+c_{2}\left(  1+h_{i}\right)  ^{2}-h_{i}%
^{2}\right)  }, \label{ingredient2}%
\end{equation}%
\begin{equation}
\mathrm{d}m_{x_{i}}\left(  0\right)  /\mathrm{d}x=\frac{-\left(  c_{2}\left(
1+h_{i}\right)  -h_{i}\right)  ^{2}}{\left(  h_{i}+c_{1}\right)  ^{2}\left(
c_{1}+c_{2}\left(  1+h_{i}\right)  ^{2}-h_{i}^{2}\right)  }.
\label{ingredient3}%
\end{equation}
Using (\ref{ingredient1}), (\ref{ingredient2}), and (\ref{ingredient3}) in
(\ref{normality}) and (\ref{normality Setting 2}), we obtain

\begin{proposition}
\label{Proposition2}For any $h_{1}>...>h_{m}>\bar{h}\equiv\left(
c_{2}+r\right)  /\left(  1-c_{2}\right)  $, the joint asymptotic distribution
of $\sqrt{p}\left(  \lambda_{pi}-x_{pi}\right)  ,$ $i=1,...,m$ is normal with
diagonal covariance matrix. For Setting 1,%
\begin{equation}
\sqrt{p}\left(  \lambda_{pi}-x_{pi}\right)  \overset{d}{\rightarrow}N\left(
0,2r^{2}\frac{h_{i}^{2}\left(  h_{i}+1\right)  ^{2}\left(  h_{i}^{2}%
-c_{2}\left(  h_{i}+1\right)  ^{2}-c_{1}\right)  }{\left(  c_{2}-h_{i}%
+c_{2}h_{i}\right)  ^{4}}\right)  , \label{convergence setting 1}%
\end{equation}
whereas for Setting 2,%
\begin{equation}
\sqrt{p}\left(  \lambda_{pi}-x_{pi}\right)  \overset{d}{\rightarrow}N\left(
0,2t^{2}\frac{h_{i}^{2}\left(  h_{i}+1\right)  ^{2}\left(  h_{i}^{2}%
-c_{2}\left(  h_{i}+1\right)  ^{2}-c_{1}\right)  }{\left(  c_{2}-h_{i}%
+c_{2}h_{i}\right)  ^{4}}\right)  . \label{convergence setting 2}%
\end{equation}
Here%
\begin{align*}
r^{2}  &  =c_{1}+c_{2}-c_{1}c_{2},\\
t^{2}  &  =c_{1}+c_{2}-\frac{c_{1}\left(  h_{i}^{2}-c_{1}\right)  }{\left(
1+h_{i}\right)  ^{2}},
\end{align*}
and%
\[
x_{pi}=\frac{\left(  h_{i}+p/n_{1}\right)  \left(  h_{i}+1\right)  }%
{h_{i}-\left(  h_{i}+1\right)  p/n_{2}}.
\]

\end{proposition}

\begin{remark}
It is straightforward to verify that $t^{2}<r^{2}$ as long as $h_{i}>\bar{h}$.
Therefore, the asymptotic variance of $\lambda_{i}$ is smaller for Setting 2
than for Setting 1. This accords with intuition because, as discussed above,
Setting 2 corresponds to the asymptotic analysis conditional on factors
$\mathcal{F},$ whereas Setting 1 corresponds to the unconditional asymptotic
analysis. The factors' variance adds to the asymptotic variance of
$\lambda_{i}$.
\end{remark}

\begin{remark}
For Setting 1, when $c_{2}\rightarrow0,$ the asymptotic variance of
$\lambda_{i}$ converges to the correct asymptotic variance%
\[
2c_{1}\left(  h_{i}+1\right)  ^{2}\left(  h_{i}^{2}-c_{1}\right)  /h_{i}^{2}%
\]
of the largest eigenvalue of the spiked Wishart model. Non-centrality spikes
in Wishart distribution were considered in Onatski (2007). The limit of the
asymptotic variance in (\ref{convergence setting 2}) when $c_{2}\rightarrow0$
coincides with the formula for the asymptotic variance derived there.
\end{remark}

\section{Analysis of the joint density of eigenvalues}

From now on, let us consider the case of a single spike, which is located
above the phase transition threshold $\bar{h}$. That is, assume that $k=m=1,$
and let $h_{1}=h_{p}$. We would like to study the asymptotic behavior of the
ratio of the joint densities of all the eigenvalues of $\mathbf{F}$ that
correspond to%
\[
H_{0}:h_{p}=h_{0}\text{ and to }H_{1}:h_{p}=h_{0}+\gamma/\sqrt{p},
\]
where $h_{0}>\bar{h}$ is fixed and $\gamma$ is a local parameter.

Following James (1964) and Khatri (1967), we can write the joint density of
the eigenvalues of $\mathbf{F}$ in Setting 1 as
\[
f_{1}(\Lambda;h_{p})=\frac{Z_{p1}(\Lambda)}{(1+h_{p})^{n_{A}/2}}\left.
_{1}F_{0}\right.  \left(  \frac{n}{2};\frac{h_{p}}{h_{p}+1}VV^{\prime}%
,\alpha_{p}\Lambda(I_{p}+\alpha_{p}\Lambda)^{-1}\right)  ,
\]
and in Setting 2 as%
\[
f_{2}(\Lambda;h)=\frac{Z_{p2}(\Lambda)}{\exp\left\{  h_{p}n_{A}/2\right\}
}\left.  _{1}F_{1}\right.  \left(  \frac{n}{2};\frac{n_{A}}{2};\frac
{n_{A}h_{p}}{2}VV^{\prime},\alpha_{p}\Lambda(I_{p}+\alpha_{p}\Lambda
)^{-1}\right)  ,
\]
where $\left.  _{1}F_{0}\right.  $ and $\left.  _{1}F_{1}\right.  $ are the
hypergeometric functions of two matrix arguments, $\alpha_{p}=n_{A}/n_{2}$,
$n=n_{A}+n_{2}$, $\Lambda=\operatorname*{diag}\left\{  \lambda_{p1}%
,\cdots,\lambda_{pp}\right\}  $, and $Z_{pj}(\Lambda),$ $j=1,2,$ depend on
$n_{A},n_{2},p$ and $\Lambda$, but not on $h_{p}$. The joint densities are
evaluated at the observed values of the eigenvalues.

To facilitate analysis, we use Proposition 1 of Dharmawansa and Johnstone
(2014) to rewrite $f_{1}(\Lambda;h_{p})$ and $f_{2}(\Lambda;h_{p})$ as shown
in the following lemma.

\begin{lemma}
\label{Lem1} Consider the region $\mathbb{C}\backslash(1,\infty)$ in the
complex plane. Let $\mathcal{\tilde{K}}$ be a contour defined in that region
which starts at $-\infty$, encircles $\tilde{\lambda}_{pj}=\alpha_{p}%
\lambda_{pj}/\left(  1+\alpha_{p}\lambda_{pj}\right)  ,$ $j=1,...,p,$
counter-clockwise and returns to $-\infty$. Then we have
\begin{equation}
f_{1}(\Lambda;h_{p})=C_{p1}(\Lambda)k_{p1}(h_{p})\frac{1}{2\pi i}%
\int_{\mathcal{\tilde{K}}}\left(  1-\frac{h_{p}}{1+h_{p}}z\right)
^{\frac{p-n-2}{2}}\prod_{j=1}^{p}\left(  z-\tilde{\lambda}_{pj}\right)
^{-\frac{1}{2}}\mathrm{d}z, \label{eig_pdf_1}%
\end{equation}
and%
\begin{align}
f_{2}(\Lambda;h_{p})  &  =C_{p2}(\Lambda)k_{p2}(h_{p})\frac{1}{2\pi
i}\label{eig_pdf_2}\\
&  \times\int_{\mathcal{\tilde{K}}}\left.  _{1}F_{1}\left(  \frac{n-p+2}%
{2},\frac{n_{A}-p+2}{2};\frac{n_{A}h_{p}}{2}z\right)  \right.  \prod_{j=1}%
^{p}\left(  z-\tilde{\lambda}_{pj}\right)  ^{-\frac{1}{2}}\mathrm{d}z\nonumber
\end{align}
where $C_{pj}(\Lambda),$ $j=1,2,$ depend on $n_{A},n_{2},p$ and $\Lambda,$ but
not on $h_{p}$, $k_{p1}(h_{p})=(1+h_{p})^{\frac{p-2-n_{A}}{2}}h_{p}^{1-p/2}$,
and $k_{p2}(h_{p})=\exp\left\{  -n_{A}h_{p}/2\right\}  h_{p}^{1-p/2}.$
\end{lemma}

We will now derive an asymptotic approximation to the contour integrals in
(\ref{eig_pdf_1}) and (\ref{eig_pdf_2}). First, we will analyze
(\ref{eig_pdf_1}) and then turn to (\ref{eig_pdf_2}).

\subsection{Asymptotic approximation: Setting 1}

\begin{figure}[ptb]%
\centering
\includegraphics[
height=3.3434in,
width=4.2021in
]%
{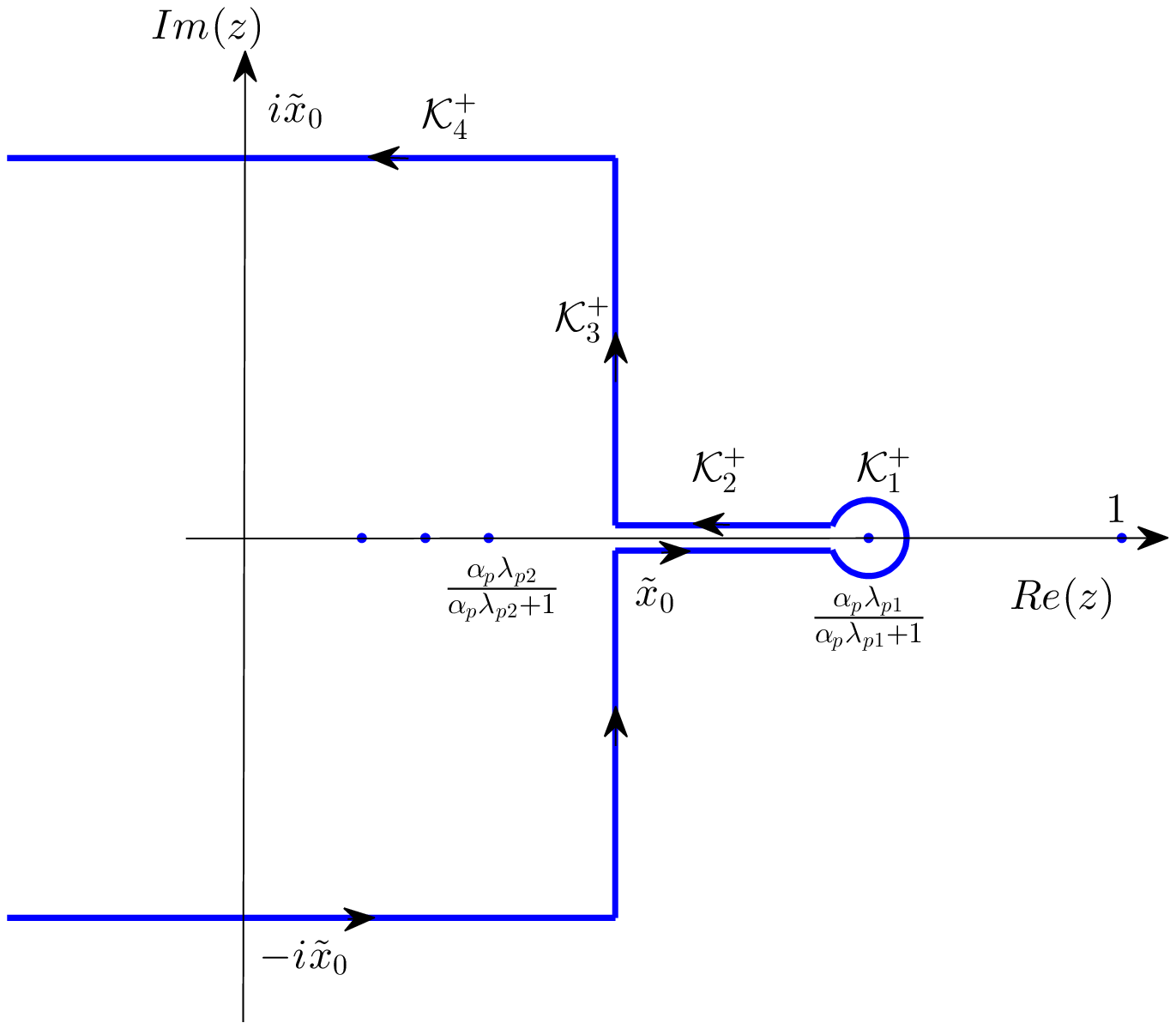}%
\caption{The contour $\mathcal{K}$.}%
\label{figure1}%
\end{figure}

Let us deform the contour $\mathcal{\tilde{K}}$, without changing the
integral's value with probability approaching one as $p,\mathbf{n}%
\rightarrow_{\mathbf{c}}\infty$, as shown in Figure \ref{figure1}. Let
$\mathcal{K}=\mathcal{K}^{+}\cup\overline{\mathcal{K}^{+}}$ with
$\mathcal{K}^{+}=\mathcal{K}_{1}^{+}\cup\mathcal{K}_{2}^{+}\cup\mathcal{K}%
_{3}^{+}\cup\mathcal{K}_{4}^{+}$, where
\begin{align*}
\mathcal{K}_{1}^{+}  &  =\left\{  z:\text{$\Im(z)\geq0$ and $|z-\tilde
{\lambda}_{p1}|=\epsilon$}\right\}  ,\\
\mathcal{K}_{2}^{+}  &  =\{z:z\in\lbrack\tilde{x}_{0},\tilde{\lambda}%
_{p1}-\epsilon]\},\\
\mathcal{K}_{3}^{+}  &  =\{z:\text{$\Re(z)=\tilde{x}_{0}$ and $0\leq\Im
(z)\leq\tilde{x}_{0}$}\},\;\;\text{and}\\
\mathcal{K}_{4}^{+}  &  =\{z:\text{$\Re(z)\leq\tilde{x}_{0}$ and
$\Im(z)=\tilde{x}_{0}$}\}.
\end{align*}
Here $\epsilon>0$ is a small number and $\alpha b_{+}/\left(  1+\alpha
b_{+}\right)  <\tilde{x}_{0}<\alpha x_{1}/\left(  1+\alpha x_{1}\right)  $
with
\[
\alpha=\lim\alpha_{p}=c_{2}/c_{1}\text{,}%
\]
and%
\begin{equation}
x_{1}=\lim x_{p1}=\frac{\left(  h_{0}+c_{1}\right)  \left(  h_{0}+1\right)
}{h_{0}-\left(  h_{0}+1\right)  c_{2}}. \label{def_x1}%
\end{equation}
As follows from our results in the previous section, $\tilde{\lambda}%
_{p1}\overset{a.s.}{\rightarrow}\alpha x_{1}/\left(  1+\alpha x_{1}\right)  $
and $\tilde{\lambda}_{p2}\overset{a.s.}{\rightarrow}\alpha b_{+}/\left(
1+\alpha b_{+}\right)  $, so $\tilde{x}_{0}\in\left(  \tilde{\lambda}%
_{p2},\tilde{\lambda}_{p1}\right)  $ for sufficiently large $p$ and
$\mathbf{n},$ a.s.%

Consider the following integral over the deformed contour $\mathcal{K}$%
\begin{equation}
I_{p}(\gamma,\Lambda)=\int_{\mathcal{K}}\left.  _{1}\mathcal{F}_{0}\left(
z\right)  \right.  \prod_{j=1}^{p}\left(  z-\tilde{\lambda}_{pj}\right)
^{-\frac{1}{2}}\mathrm{d}z, \label{Ip_def}%
\end{equation}
where
\[
\left.  _{1}\mathcal{F}_{0}\left(  z\right)  \right.  \equiv\left(
1-\frac{h_{p}}{1+h_{p}}z\right)  ^{\frac{p-n-2}{2}}.
\]
For two sequences of random variables $\left\{  \xi_{p}\right\}  $ and
$\left\{  \eta_{p}\right\}  ,$ we will write $\xi_{p}\overset{\mathrm{P}%
}{\sim}\eta_{p}$ if and only if $\xi_{p}/\eta_{p}$ converges in probability to
1 as $p,\mathbf{n}\rightarrow_{\mathbf{c}}\infty$. We have the following lemma.

\begin{lemma}
\label{Lem2}Under the hypothesis that $h_{p}=h_{0}$, uniformly in $\gamma$
from any compact subset of $\mathbb{R}$
\[
I_{p}(\gamma,\Lambda)\overset{\mathrm{P}}{\sim}\left.  _{1}\mathcal{F}%
_{0}\left(  \tilde{\lambda}_{p1}\right)  \right.  2\left(  -\frac{\pi}{pH_{0}%
}\right)  ^{\frac{1}{2}}\prod_{j=2}^{p}\left(  \tilde{\lambda}_{p1}%
-\tilde{\lambda}_{pj}\right)  ^{-\frac{1}{2}},
\]
where the principal branches of the square roots are used, and
\begin{equation}
H_{0}=\frac{h_{0}(1-c_{2})(\mu_{0}-\sqrt{b_{+}})(\mu_{0}-\sqrt{b_{-}}%
)(c_{1}+c_{2}\mu_{0})}{2c_{1}c_{2}\left(  h_{0}-c_{2}\mu_{0}\right)  \mu
_{0}(c_{1}+h_{0})} \label{def_H0}%
\end{equation}
with $\mu_{0}=h_{0}+1$.
\end{lemma}

\textbf{Proof:} Let $\mathcal{K}_{j}=\mathcal{K}_{j}^{+}\cup\overline
{\mathcal{K}_{j}^{+}}$ for $j=1,...,4.$ Using this notation, we can decompose
(\ref{Ip_def}) as
\begin{equation}
I_{p}(\gamma,\Lambda)=I_{1,2,p}(\gamma,\Lambda)+I_{3,4,p}(\gamma,\Lambda)
\label{intdecom}%
\end{equation}
where $I_{1,2,p}(\gamma,\Lambda)$ is the part of the integral corresponding to
$\mathcal{K}_{1}\cup\mathcal{K}_{2},$ and $I_{3,4,p}(\theta,\lambda)$ is the
part corresponding to the rest of the contour, $\mathcal{K}_{3}\cup
\mathcal{K}_{4}$. Our strategy is to show that the integral $I_{p}%
(\gamma,\Lambda)$ is asymptotically equivalent to $I_{1,2,p}(\gamma,\Lambda)$,
the integral $I_{3,4,p}(\gamma,\Lambda)$ being asymptotically dominated by
$I_{1,2,p}(\gamma,\Lambda)$.

Let us first focus on $I_{1,2,p}(\gamma,\Lambda)$. Since the singularity of
the integrand at $\tilde{\lambda}_{p1}$ is of the inverse square root type, as
the radius $\epsilon$ of $\mathcal{K}_{1}$ converges to zero, the integral
over $\mathcal{K}_{1}$ converges to zero too. Therefore, we have%
\begin{equation}
I_{1,2,p}(\gamma,\Lambda)=2I_{2,p}\left(  \gamma,\Lambda\right)  ,
\label{cont1decom}%
\end{equation}
where%
\[
I_{2,p}\left(  \gamma,\Lambda\right)  =\int_{\tilde{\lambda}_{p1}}^{\tilde
{x}_{0}}\left.  _{1}\mathcal{F}_{0}\left(  z\right)  \right.  \prod_{j=1}%
^{p}\left(  z-\tilde{\lambda}_{pj}\right)  ^{-\frac{1}{2}}\mathrm{d}z.
\]

Changing the variable of integration from $z$ to $x=\tilde{\lambda}_{p1}-z,$
we arrive at
\begin{align*}
I_{2,p}(\theta,\lambda)  &  =-\int_{0}^{\tilde{\lambda}_{p1}-\tilde{x}_{0}%
}\left.  _{1}\mathcal{F}_{0}\left(  \tilde{\lambda}_{p1}-x\right)  \right.
\prod_{j=1}^{p}\left(  \tilde{\lambda}_{p1}-\tilde{\lambda}_{pj}-x\right)
^{-\frac{1}{2}}\mathrm{d}x\\
&  =\left.  _{1}\mathcal{F}_{0}\left(  \tilde{\lambda}_{p1}\right)  \right.
(-1)^{\frac{1}{2}}\int_{0}^{\tilde{\lambda}_{p1}-\tilde{x}_{0}}\exp\left\{
-pf_{p}(x)\right\}  x^{-\frac{1}{2}}\mathrm{d}x
\end{align*}
where
\[
f_{p}(x)=\frac{\beta_{p}}{2}\ln\left(  1+\frac{h_{p}(1+\alpha_{p}\lambda
_{p1})}{1+h_{p}+\alpha_{p}\lambda_{p1}}x\right)  +\frac{1}{2p}\sum_{j=2}%
^{p}\ln\left(  \tilde{\lambda}_{p1}-\tilde{\lambda}_{pj}-x\right)
\]
with $\beta_{p}=\left(  2+n-p\right)  /p$. Now the integral $\int_{0}%
^{\tilde{\lambda}_{p1}-\tilde{x}_{0}}\exp\left\{  -pf_{p}(x)\right\}
x^{-\frac{1}{2}}\mathrm{d}x$ can be evaluated using standard Laplace
approximation steps (see Olver (1997), section 7.3) as follows.

First, let us show that the derivative $\frac{\mathrm{d}}{\mathrm{d}x}%
f_{p}(x)$ is continuous and positive on $x\in\lbrack0,\tilde{\lambda}%
_{p1}-\tilde{x}_{0}]$ for sufficiently large $p$ and $\mathbf{n},$ a.s. We
have
\[
\frac{\mathrm{d}}{\mathrm{d}x}f_{p}(x)=\frac{\beta_{p}}{2}\frac{h_{p}\left(
1+\alpha_{p}\lambda_{p1}\right)  }{\left(  h_{p}+(1+\alpha_{p}\lambda
_{p1})(1+h_{p}x)\right)  }-\frac{1}{2p}\sum_{j=2}^{p}\left(  \tilde{\lambda
}_{p1}-\tilde{\lambda}_{pj}-x\right)  ^{-1}.
\]
Therefore, the continuity follows from the fact that, when $h_{p}=h_{0}$,
\[
\lambda_{p1}\overset{a.s.}{\rightarrow}x_{1}\equiv\frac{\left(  h_{0}%
+c_{1}\right)  \left(  h_{0}+1\right)  }{h_{0}-\left(  h_{0}+1\right)  c_{2}%
}\text{ and }\lambda_{p2}\overset{a.s.}{\rightarrow}b_{+}.
\]
In order to establish the positivity, we first obtain
\[
\min_{x\in\lbrack0,\tilde{\lambda}_{p1}-\tilde{x}_{0}]}\frac{\mathrm{d}%
}{\mathrm{d}x}f_{p}(x)=\frac{\beta_{p}}{2}\frac{h_{p}}{1+h_{p}-h_{p}\tilde
{x}_{0}}-\frac{1}{2p}\sum_{j=2}^{p}\left(  \tilde{x}_{0}-\tilde{\lambda}%
_{pj}\right)  ^{-1}.
\]
It is straightforward to verify that the above equation can be represented in
the following form
\begin{align*}
\min_{x\in\lbrack0,\tilde{\lambda}_{p1}-\tilde{x}_{0}]}\frac{\mathrm{d}%
}{\mathrm{d}x}f_{p}(x)  &  =\frac{\beta_{p}}{2}\frac{h_{p}(1+\alpha x_{0}%
)}{1+h_{p}+\alpha x_{0}}+\frac{p-1}{2p}(1+\alpha x_{0})\\
&  +\frac{1}{2\alpha_{p}p}(1+\alpha x_{0})^{2}\sum_{j=2}^{p}\left(
\lambda_{pj}-\alpha x_{0}/\alpha_{p}\right)  ^{-1}.
\end{align*}
where $x_{0}=\tilde{x}_{0}/\left(  \alpha(1-\tilde{x}_{0})\right)  .$
Therefore, we obtain%
\begin{equation}
\min_{x\in\lbrack0,\tilde{\lambda}_{p1}-\tilde{x}_{0}]}\frac{\mathrm{d}%
}{\mathrm{d}x}f_{p}(x)\overset{a.s.}{\rightarrow}\Psi(x_{0},h_{0}) \label{as1}%
\end{equation}
where
\[
\Psi(x_{0},h_{0})=\frac{\beta h_{0}(1+\alpha x_{0})}{2(1+h_{0}+\alpha x_{0}%
)}+\frac{1}{2}(1+\alpha x_{0})+\frac{1}{2\alpha}(1+\alpha x_{0})^{2}m\left(
x_{0}\right)  ,
\]%
\[
\beta=c_{1}^{-1}+c_{2}^{-1}-1,
\]
and $m(x_{0})=\lim_{z\rightarrow x_{0}}m\left(  z\right)  $ with $m(z)$ being
the Stieltjes transform of the limiting spectral distribution of $\mathbf{F},$
that is the distribution with density (\ref{Wachter density}).

Since $m(x_{0})$ is increasing on $x_{0}\in(b_{+},\infty)$, we have
\[
\Psi(x_{0},h_{0})>\lim_{x_{0}\downarrow b_{+}}\Psi(x_{0},h_{0}).
\]
Moreover, noting the fact that $\lim_{x_{0}\downarrow b_{+}}\Psi(x_{0},h_{0})$
is an increasing function of $h_{0}$ and $h_{0}>\bar{h}\equiv\left(
c_{2}+r\right)  /\left(  1-c_{2}\right)  =\sqrt{b_{+}}-1$, we obtain
\begin{equation}
\Psi(x_{0},h_{0})>\lim_{x_{0}\downarrow b_{+}}\Psi(x_{0},h_{0})>\lim
_{x_{0}\downarrow b_{+}}\Psi\left(  x_{0},\bar{h}\right)  . \label{as2}%
\end{equation}
Finally, direct calculations, which are not reported here to save space, show
that, as $x_{0}$ converges to $b_{+}$ from the right,
\begin{equation}
m(x_{0})\rightarrow-1/(b_{+}-\sqrt{b_{+}}). \label{Prathapas result}%
\end{equation}

This in turn gives
\begin{equation}
\lim_{x\downarrow b_{+}}\Psi\left(  x,\bar{h}\right)  =0 \label{as3}%
\end{equation}
which establishes the positivity.

Since $\lambda_{p1}\overset{a.s.}{\rightarrow}x_{1}$, we have
\[
f_{p}^{^{\prime}}(0)\equiv\left.  \frac{\mathrm{d}}{\mathrm{d}x}%
f_{p}(x)\right\vert _{x=0}\overset{a.s.}{\rightarrow}H_{0},
\]
where%
\[
H_{0}=\frac{r^{2}}{2c_{1}c_{2}}\frac{(1+\alpha x_{1})h_{0}}{1+h_{0}+\alpha
x_{1}}+\frac{1}{2}(1+\alpha x_{1})+\frac{1}{2\alpha}(1+\alpha x_{1}%
)^{2}m(x_{1}).
\]
Direct calculations show that
\begin{equation}
m(x_{1})=\lim_{z\rightarrow x_{1}}m(z)=-(1+h_{0})/\left(  x_{1}h_{0}\right)  ,
\label{Iains result}%
\end{equation}
which, after some algebraic manipulations, gives (\ref{def_H0}).

We may now exploit the approach given in Olver (1997, pp. 81-82) to yield
\[
\int_{0}^{\tilde{\lambda}_{p1}-\tilde{x}_{0}}e^{-pf_{p}(x)}x^{-\frac{1}{2}%
}\mathrm{d}x\overset{\mathrm{P}}{\sim}\left(  \frac{\pi}{pf_{p}^{^{\prime}%
}(0)}\right)  ^{\frac{1}{2}}e^{-pf_{p}(0)}.
\]
Therefore, we obtain
\begin{equation}
I_{2,p}(\gamma,\Lambda)\overset{\mathrm{P}}{\sim}\left.  _{1}\mathcal{F}%
_{0}\left(  \tilde{\lambda}_{p1}\right)  \right.  \left(  -\frac{\pi}{pH_{0}%
}\right)  ^{\frac{1}{2}}\prod_{j=2}^{p}\left(  \tilde{\lambda}_{p1}%
-\tilde{\lambda}_{pj}\right)  ^{-\frac{1}{2}}. \label{cont1ans}%
\end{equation}
As Lemma \ref{Lem:negl} below shows, $I_{3,4,p}(\gamma,\Lambda)$ is
asymptotically dominated by $I_{2,p}(\theta,\lambda)$, which completes the
proof.$\square$

\begin{lemma}
\label{Lem:negl} Under the hypothesis that $h_{p}=h_{0},$ uniformly in
$\gamma$ from any compact subset of $\mathbb{R}$
\begin{equation}
I_{3,4,p}(\gamma,\Lambda)=o_{\mathrm{P}}\left(  I_{2,p}(\gamma,\Lambda
)\right)  . \label{contdim}%
\end{equation}

\end{lemma}

\textbf{Proof:} Let us first consider the integral over the contour
$\mathcal{K}_{3}$. For $z\in\mathcal{K}_{3}$, we have
\[
\left\vert \left.  _{1}\mathcal{F}_{0}\left(  z\right)  \right.  \prod
_{j=1}^{p}\left(  z-\tilde{\lambda}_{pj}\right)  ^{-\frac{1}{2}}\right\vert
<\left.  _{1}\mathcal{F}_{0}\left(  \tilde{\lambda}_{p1}\right)  \right.
e^{-pf_{p}\left(  \tilde{\lambda}_{p1}-\tilde{x}_{0}\right)  }\left(
\tilde{\lambda}_{p1}-\tilde{x}_{0}\right)  ^{-\frac{1}{2}}.
\]
Also, in view of (\ref{as1}), (\ref{as2}), and (\ref{as3}), we have
$f_{p}(\tilde{\lambda}_{1p}-\tilde{x}_{0})>f_{p}(0)+\epsilon$, for
sufficiently large $p$ and $\mathbf{n}$, a.s., where $\epsilon>0$. Therefore,
using (\ref{cont1ans}), we conclude
\begin{equation}
\int_{\mathcal{K}_{3}}\left.  _{1}\mathcal{F}_{0}\left(  z\right)  \right.
\prod_{j=1}^{p}\left(  z-\tilde{\lambda}_{pj}\right)  ^{-\frac{1}{2}%
}\mathrm{d}z=o_{\mathrm{P}}\left(  I_{2,p}(\gamma,\Lambda)\right)  .
\label{cont3}%
\end{equation}

Now consider the integral over the contour $\mathcal{K}_{4}$. We have
\begin{align*}
&  \left\vert \int_{\mathcal{K}_{4}}\left.  _{1}\mathcal{F}_{0}\left(
z\right)  \right.  \prod_{j=1}^{p}\left(  z-\tilde{\lambda}_{pj}\right)
^{-\frac{1}{2}}\mathrm{d}z\right\vert <2\prod_{j=1}^{p}\left\vert
\tilde{\lambda}_{pj}-\tilde{x}_{0}\right\vert ^{-\frac{1}{2}}\int_{-\infty
}^{\tilde{x}_{0}}\left.  _{1}\mathcal{F}_{0}\left(  x\right)  \right.
\mathrm{d}x\\
&  \qquad\qquad=4\frac{1+h_{p}}{h_{p}(n-p)}\left(  1-\frac{h_{p}}{1+h_{p}%
}\tilde{x}_{0}\right)  \left.  _{1}\mathcal{F}_{0}\left(  \tilde{x}%
_{0}\right)  \right.  \prod_{j=1}^{p}\left\vert \tilde{\lambda}_{pj}-\tilde
{x}_{0}\right\vert ^{-\frac{1}{2}}.
\end{align*}
Since%
\[
\left.  _{1}\mathcal{F}_{0}\left(  \tilde{x}_{0}\right)  \right.  \prod
_{j=1}^{p}\left(  \tilde{\lambda}_{pj}-\tilde{x}_{0}\right)  ^{-\frac{1}{2}%
}=\left.  _{1}\mathcal{F}_{0}\left(  \tilde{\lambda}_{p1}\right)  \right.
e^{-pf_{p}\left(  \tilde{\lambda}_{p1}-\tilde{x}_{0}\right)  }\left(
\tilde{\lambda}_{p1}-\tilde{x}_{0}\right)  ^{-\frac{1}{2}},
\]
we can follow a similar procedure to that outlined above to obtain
\[
\int_{\mathcal{K}_{4}}\left.  _{1}\mathcal{F}_{0}\left(  z\right)  \right.
\prod_{j=1}^{p}\left(  z-\tilde{\lambda}_{pj}\right)  ^{-\frac{1}{2}%
}\mathrm{d}z=o_{\mathrm{P}}\left(  I_{2,p}(\gamma,\Lambda)\right)  .
\]
This along with (\ref{cont3}) gives (\ref{contdim}).$\square$

\subsection{Asymptotic approximation: Setting 2}

Consider the following integral%
\[
J_{p}(\gamma,\Lambda)=\int_{\mathcal{\tilde{K}}}\left.  _{1}\mathcal{F}%
_{1}\left(  \zeta\right)  \right.  \prod_{j=1}^{p}\left(  z-\tilde{\lambda
}_{pj}\right)  ^{-\frac{1}{2}}\mathrm{d}z,
\]
where%
\[
\left.  _{1}\mathcal{F}_{1}\left(  \zeta\right)  \right.  \equiv\left.
_{1}F_{1}\left(  n_{A}u+1,n_{A}v+1;n_{A}\zeta\right)  \right.
\]
with
\[
u=\frac{n_{A}+n_{2}-p}{2n_{A}},\text{ }v=\frac{n_{A}-p}{2n_{A}},\text{ and
}\zeta=\frac{h_{p}}{2}z.
\]
In Johnstone and Onatski (2014) (Theorem 5), the following result is derived.
As \thinspace$p,\mathbf{n\rightarrow}_{\mathbf{c}}\infty$,%

\begin{equation}
\left.  _{1}\mathcal{F}_{1}\left(  \zeta\right)  \right.  =\frac{C\left(
p,\mathbf{n}\right)  e^{-n_{A}\varphi\left(  \zeta\right)  }}{\sqrt{2\pi
n_{A}}i}\left(  \psi\left(  \zeta\right)  +O\left(  n_{A}^{-1}\right)
\right)  , \label{approx F11}%
\end{equation}
where $O\left(  n_{A}^{-1}\right)  $ is uniform for $\zeta$ that do not
approach zero or negative semi-axis and $\mathfrak{R}\left(  \zeta\right)
\geq-2u+v$,%
\[
C\left(  p,\mathbf{n}\right)  =\frac{\Gamma\left(  n_{A}v+1\right)
\Gamma\left(  n_{A}\left(  u-v\right)  +1\right)  }{\Gamma\left(
n_{A}u+1\right)  },
\]%
\begin{equation}
\varphi\left(  \zeta\right)  =\left(  u-v\right)  \ln\left(  u-v\right)
+v\ln\left(  z_{+}+v\right)  -u\ln\left(  z_{+}+u\right)  -z_{+},
\label{def_fi}%
\end{equation}
where the principal branches of the logarithms are chosen,%
\begin{equation}
z_{+}=\frac{1}{2}\left\{  \zeta-v+\sqrt{\left(  \zeta-v\right)  ^{2}+4u\zeta
}\right\}  , \label{def z+}%
\end{equation}
where the principal branch of the square root is chosen when $\mathfrak{R}%
\left(  \zeta\right)  \geq-2u+v$ and the other branch is chosen when
$\operatorname{Re}\zeta<-2u+v$, and%
\[
\psi\left(  \zeta\right)  =\left[  \left(  z_{+}-\zeta\right)  \sqrt{\frac
{u}{z_{+}^{2}}-\frac{u-v}{\left(  z_{+}-\zeta\right)  ^{2}}}\right]
^{-1},\text{ }%
\]
where the branch of the square root is chosen so that $\sqrt{-1}=-i.$

\begin{figure}[ptb]%
\centering
\includegraphics[
height=3.3434in,
width=4.2021in
]%
{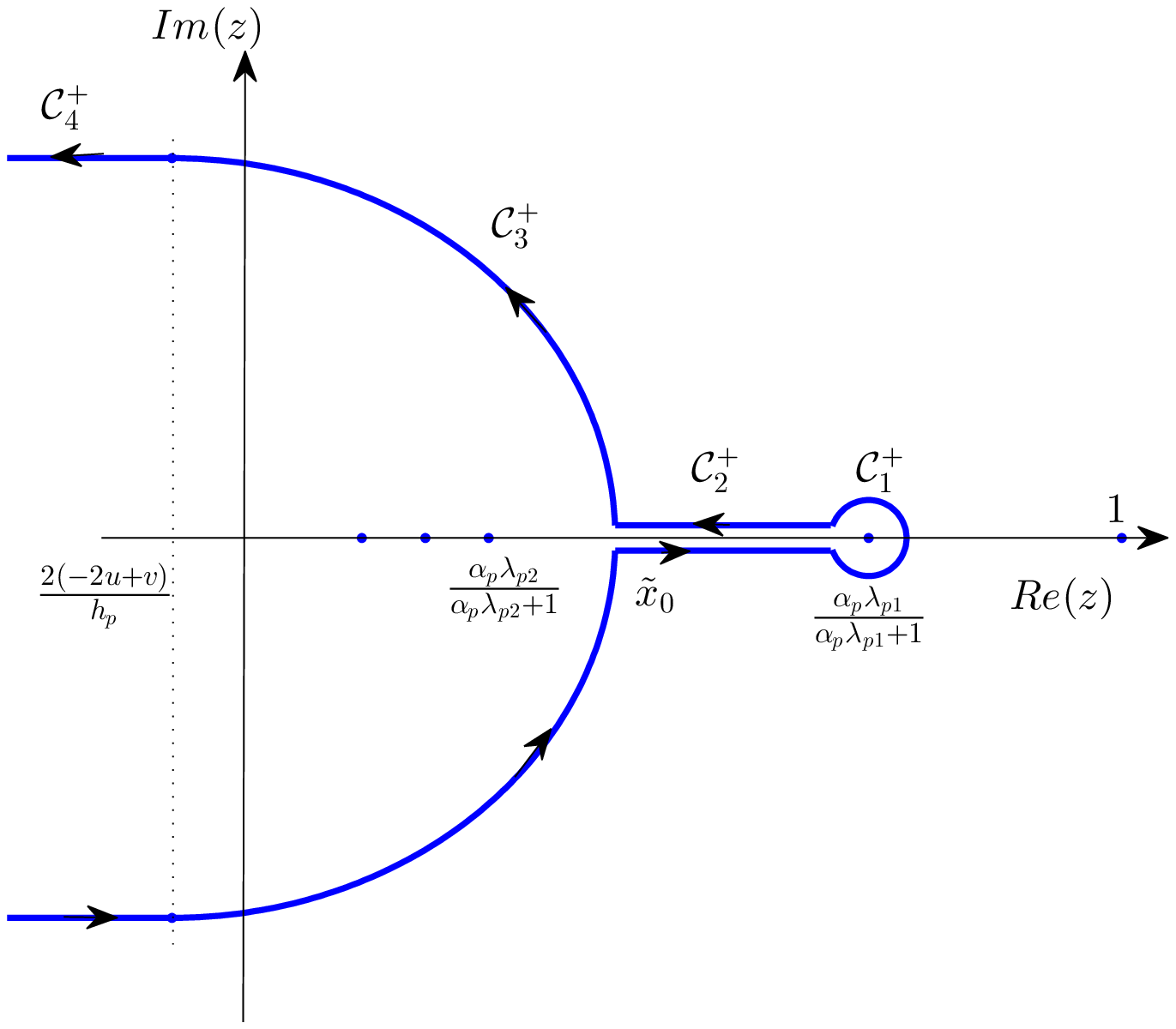}%
\caption{The contour $\mathcal{C}$.}%
\label{figure2}%
\end{figure}

We will deform the contour $\mathcal{\tilde{K}}$, without changing the
integral's value with probability approaching one as $p,\mathbf{n}%
\rightarrow_{\mathbf{c}}\infty$, as shown in Figure \ref{figure2}. Formally,
$\mathcal{C}=\mathcal{C}^{+}\cup\overline{\mathcal{C}^{+}}$ with
$\mathcal{C}^{+}=\mathcal{C}_{1}^{+}\cup\mathcal{C}_{2}^{+}\cup\mathcal{C}%
_{3}^{+}\cup\mathcal{C}_{4}^{+}$, where
\begin{align*}
\mathcal{C}_{1}^{+}  &  =\{z:\text{$\Im(z)\geq0$ and $|z-\tilde{\lambda}%
_{p1}|=\epsilon\}$},\\
\mathcal{C}_{2}^{+}  &  =\{z:z\in\lbrack\tilde{x}_{0},\tilde{\lambda}%
_{p1}-\epsilon]\},\\
\mathcal{C}_{3}^{+}  &  =\left\{  z:z=2\zeta/h_{p}\text{ s.t. }\mathfrak{R}%
\left(  \zeta\right)  \geq-2u+v,\text{ $\Im$}\left(  \zeta\right)
\geq0,\text{ }\left\vert z_{+}+u\right\vert =\left\vert z_{0+}+u\right\vert
\right\}  ,\\
\mathcal{C}_{4}^{+}  &  =\{z:z=2\zeta/h_{p}\text{ s.t. }\mathfrak{R}\left(
\zeta\right)  <-2u+v,\text{ and $\Im$}\left(  z_{+}\right)  =\left\vert
z_{0+}+u\right\vert \}.
\end{align*}
Here%
\begin{align}
z_{0+}  &  =\frac{1}{2}\left\{  \zeta_{0}-v+\sqrt{\left(  \zeta_{0}-v\right)
^{2}+4u\zeta_{0}}\right\}  ,\text{ and}\label{def z0+}\\
\zeta_{0}  &  =\frac{h_{p}}{2}\tilde{x}_{0}.\nonumber
\end{align}
%

\begin{lemma}
\label{Lem2a}Under the hypothesis that $h_{p}=h_{0}$, uniformly in $\gamma$
from any compact subset of $\mathbb{R}$%
\[
J_{p}(\gamma,\Lambda)\overset{\mathrm{P}}{\sim}2Z\left(  p,\mathbf{n}%
,h_{0}\right)  e^{-n_{A}\varphi\left(  \frac{h_{p}}{2}\tilde{\lambda}%
_{p1}\right)  }\left(  -\frac{\pi}{pH_{0}}\right)  ^{1/2}\prod_{j=2}%
^{p}\left(  \tilde{\lambda}_{p1}-\tilde{\lambda}_{pj}\right)  ^{-\frac{1}{2}%
},
\]
where%
\[
Z\left(  p,\mathbf{n},h_{0}\right)  =\frac{C\left(  p,\mathbf{n}\right)
}{\sqrt{\pi p}}\frac{c_{1}+c_{2}\mu_{0}}{\sqrt{c_{2}\mu_{0}^{2}+c_{1}%
^{2}-c_{1}+2c_{1}\mu_{0}}},
\]
and $H_{0}$ and $\mu_{0}$ are as defined in Lemma \ref{Lem2}.
\end{lemma}

\textbf{Proof}: Similar to the case of Setting 1, we split $J_{p}%
(\gamma,\Lambda)$ into two parts%
\[
J_{p}(\gamma,\Lambda)=J_{1,2,p}(\gamma,\Lambda)+J_{3,4,p}(\gamma,\Lambda),
\]
where $J_{1,2,p}(\gamma,\Lambda)$ is the part of the integral corresponding to
$\mathcal{C}_{1}\cup\mathcal{C}_{2},$ and $J_{3,4,p}(\theta,\lambda)$ is the
part corresponding to the rest of the contour, $\mathcal{C}_{3}\cup
\mathcal{C}_{4}$. Furthermore,%
\begin{equation}
J_{1,2,p}(\gamma,\Lambda)\overset{\mathrm{P}}{\sim}2J_{2,p}\left(
\gamma,\Lambda\right)  , \label{cont1decom2}%
\end{equation}
where%
\[
J_{2,p}\left(  \gamma,\Lambda\right)  =\int_{\tilde{\lambda}_{p1}}^{\tilde
{x}_{0}}\frac{C\left(  p,\mathbf{n}\right)  e^{-n_{A}\varphi\left(
\zeta\right)  }}{\sqrt{2\pi n_{A}}i}\psi\left(  \zeta\right)  \prod_{j=1}%
^{p}\left(  z-\tilde{\lambda}_{pj}\right)  ^{-\frac{1}{2}}\mathrm{d}z.
\]
In contrast to (\ref{cont1decom}), we only have the asymptotic equivalence in
(\ref{cont1decom2}) because we are using the uniform asymptotic approximation
(\ref{approx F11}) to define $J_{2,p}\left(  \gamma,\Lambda\right)  .$

After the change of the variable of integration, $\zeta\mapsto x=\frac{h_{p}%
}{2}\tilde{\lambda}_{p1}-\zeta,$ we obtain%
\[
J_{2,p}\left(  \gamma,\Lambda\right)  =-\int_{0}^{\frac{h_{p}}{2}\left(
\tilde{\lambda}_{1p}-\tilde{x}_{0}\right)  }\frac{C\left(  p,\mathbf{n}%
\right)  e^{-n_{A}\varphi\left(  \zeta\right)  }}{\sqrt{2\pi n_{A}}i}%
\psi\left(  \zeta\right)  \frac{2}{h_{p}}\prod_{j=1}^{p}\left(  \frac{2}%
{h_{p}}\zeta-\tilde{\lambda}_{pj}\right)  ^{-\frac{1}{2}}\mathrm{d}x.
\]
This can be rewritten as%
\[
J_{2,p}\left(  \gamma,\Lambda\right)  =\left(  -\frac{2}{h_{p}}\right)
^{\frac{1}{2}}\frac{C\left(  p,\mathbf{n}\right)  }{\sqrt{2\pi n_{A}}i}%
\int_{0}^{\frac{h_{p}}{2}\left(  \tilde{\lambda}_{p1}-\tilde{x}_{0}\right)
}e^{-n_{A}g_{p}(\zeta)}\psi\left(  \zeta\right)  x^{-\frac{1}{2}}\mathrm{d}x,
\]
where%
\[
g_{p}(\zeta)=\varphi\left(  \zeta\right)  +\frac{1}{2n_{A}}%
{\displaystyle\sum\nolimits_{j=2}^{p}}
\ln\left(  \frac{2}{h_{p}}\zeta-\tilde{\lambda}_{pj}\right)  ,
\]
and%
\[
\zeta=\frac{h_{p}}{2}\tilde{\lambda}_{p1}-x.
\]

Following the approach in the above analysis in the case of Setting 1, we now
would like to show that the derivative $\frac{\mathrm{d}}{\mathrm{d}x}%
g_{p}(\frac{h_{p}}{2}\tilde{\lambda}_{p1}-x)$ is continuous and positive on
$x\in\lbrack0,\frac{h_{p}}{2}\left(  \tilde{\lambda}_{p1}-\tilde{x}%
_{0}\right)  ]$ for sufficiently large $p$ and $\mathbf{n},$ a.s. This is
equivalent to showing that $\frac{\mathrm{d}}{\mathrm{d}\zeta}g_{p}(\zeta)$ is
continuous and \textit{negative} on $\zeta\in\lbrack\frac{h_{p}}{2}\tilde
{x}_{0},\frac{h_{p}}{2}\tilde{\lambda}_{p1}]$ for sufficiently large $p$ and
$\mathbf{n},$ a.s.

It is straightforward to verify that $z_{+}$ satisfies the quadratic equation%
\[
z_{+}^{2}+(v-\zeta)z_{+}-u\zeta=0,
\]
and%
\begin{equation}
\zeta=\frac{z_{+}\left(  z_{+}+v\right)  }{z_{+}+u}. \label{zeta}%
\end{equation}
From this, and the definition (\ref{def z+}) of $z_{+},$ we obtain that,
$z_{+}>\zeta$ for positive $\zeta,$ and
\begin{equation}
\frac{\mathrm{d}}{\mathrm{d}\zeta}z_{+}=\frac{z_{+}+u}{2z_{+}+v-\zeta}%
=\frac{\left(  u+z_{+}\right)  ^{2}}{uv+2uz_{+}+z_{+}^{2}}>0.
\label{derivat of z}%
\end{equation}
On the other hand,%
\[
\frac{\mathrm{d}}{\mathrm{d}z_{+}}\varphi\left(  \zeta\right)  =\frac{v}%
{z_{+}+v}-\frac{u}{z_{+}+u}-1=-\frac{uv+2uz_{+}+z_{+}^{2}}{\left(
v+z_{+}\right)  \left(  u+z_{+}\right)  }<0
\]
Thus, $\varphi\left(  \zeta\right)  $ is strictly decreasing function of
$\zeta.$ Furthermore, it is a convex function of $\zeta>0$. Indeed,
\[
\frac{\mathrm{d}^{2}}{\mathrm{d}z_{+}^{2}}\varphi\left(  \zeta\right)
=-\frac{v}{\left(  z_{+}+v\right)  ^{2}}+\frac{u}{\left(  z_{+}+u\right)
^{2}}=\frac{\left(  u-v\right)  \left(  z_{+}^{2}-uv\right)  }{\left(
z_{+}+v\right)  ^{2}\left(  z_{+}+u\right)  ^{2}},
\]
and, using (\ref{zeta}) and (\ref{derivat of z}), we also have%
\[
\frac{\mathrm{d}^{2}}{\mathrm{d}\zeta^{2}}z_{+}=-2u\left(  u+z_{+}\right)
^{3}\frac{u-v}{\left(  uv+2uz_{+}+z_{+}^{2}\right)  ^{3}}.
\]
Therefore, we obtain%
\begin{align*}
\frac{\mathrm{d}^{2}}{\mathrm{d}\zeta^{2}}\varphi\left(  \zeta\right)   &
=\frac{\mathrm{d}^{2}}{\mathrm{d}z_{+}^{2}}\varphi\left(  \zeta\right)
\left(  \frac{\mathrm{d}}{\mathrm{d}\zeta}z_{+}\right)  ^{2}+\frac{\mathrm{d}%
}{\mathrm{d}z_{+}}\varphi\left(  \zeta\right)  \frac{\mathrm{d}^{2}%
}{\mathrm{d}\zeta^{2}}z_{+}\\
&  =\frac{\left(  u+z_{+}\right)  ^{2}\left(  u-v\right)  }{\left(
v+z_{+}\right)  ^{2}\left(  uv+2uz_{+}+z_{+}^{2}\right)  }>0.
\end{align*}
Therefore, $\varphi\left(  \zeta\right)  $ is, indeed, convex for positive
$\zeta,$ and has a continuous derivative.

Further, it is straightforward to see that
\[
w\left(  \zeta\right)  \equiv\frac{1}{2n_{A}}%
{\displaystyle\sum\nolimits_{j=2}^{p}}
\ln\left(  \frac{2}{h_{p}}\zeta-\tilde{\lambda}_{pj}\right)
\]
is a strictly increasing concave function of $\zeta>\frac{h_{p}}{2}%
\tilde{\lambda}_{p2}$. This implies that%
\begin{align*}
\max_{\zeta\in\left[  \frac{h_{p}}{2}\tilde{x}_{0},\frac{h_{p}}{2}%
\tilde{\lambda}_{p1}\right]  }\frac{\mathrm{d}}{\mathrm{d}\zeta}g_{p}\left(
\zeta\right)   &  <\left.  \frac{\mathrm{d}}{\mathrm{d}\zeta}\varphi\left(
\zeta\right)  \right\vert _{\zeta=\frac{h_{p}}{2}\tilde{\lambda}_{1p}}+\left.
\frac{\mathrm{d}}{\mathrm{d}\zeta}w\left(  \zeta\right)  \right\vert
_{\zeta=\frac{h_{p}}{2}\tilde{x}_{0}}\\
&  =\left.  -\frac{uv+2uz_{+}+z_{+}^{2}}{\left(  v+z_{+}\right)  \left(
u+z_{+}\right)  }\left(  1+\frac{u(u-v)}{2z_{+}\left(  u+z_{+}\right)
}\right)  \right\vert _{\zeta=\frac{h_{p}}{2}\tilde{\lambda}_{1p}}\\
&  -\frac{2p}{h_{p}n_{A}}\frac{p-1}{2p}(1+\alpha x_{0})\\
&  -\frac{2p}{h_{p}n_{A}}\frac{1}{2\alpha_{p}p}(1+\alpha x_{0})^{2}\sum
_{j=2}^{p}\left(  \lambda_{j}-\alpha x_{0}/\alpha_{p}\right)  ^{-1}%
\end{align*}
The right hand side of the latter equality a.s. converges to%
\[
\Pi(x_{0},h_{0})=-\frac{c_{1}}{c_{2}\left(  h_{0}+1\right)  }-1-\frac{c_{1}%
}{h_{0}}(1+\alpha x_{0})-\frac{c_{1}}{h_{0}\alpha}(1+\alpha x_{0})^{2}m\left(
x_{0}\right)  ,
\]
where $m(z)$ is the Stieltjes transform of the limiting spectral distribution
of $\mathbf{F}.$ Since $m(x_{0})$ is an increasing function of $x_{0}>b_{+}$,%
\[
\Pi(x_{0},h_{0})<\lim_{x_{0}\downarrow b_{+}}\Pi(x_{0},h_{0}).
\]
On the other hand, using (\ref{Prathapas result}), we get%
\begin{equation}
\lim_{x_{0}\downarrow b_{+}}\Pi(x_{0},h_{0})=-\frac{c_{1}}{c_{2}\left(
h_{0}+1\right)  }-1+r\frac{\left(  r+c_{2}\right)  ^{2}}{c_{2}h_{0}\left(
1-c_{2}\right)  \left(  r+1\right)  }. \label{lim of Psi}%
\end{equation}

Note that, considered as a function of $h_{0}>\bar{h},$ $\lim_{x_{0}\downarrow
b_{+}}\Pi(x_{0},h_{0})$ may have positive derivative only when $\lim
_{x_{0}\downarrow b_{+}}\Pi(x_{0},h_{0})<0.$ Indeed,%
\begin{align*}
\frac{\mathrm{d}}{\mathrm{d}h_{0}}\lim_{x_{0}\downarrow b_{+}}\Pi(x_{0}%
,h_{0})  &  =\frac{c_{1}}{c_{2}\left(  h_{0}+1\right)  ^{2}}-r\frac{\left(
r+c_{2}\right)  ^{2}}{c_{2}h_{0}^{2}\left(  1-c_{2}\right)  \left(
r+1\right)  }\\
&  <\frac{1}{h_{0}}\left(  \frac{c_{1}}{c_{2}\left(  h_{0}+1\right)  }%
-r\frac{\left(  r+c_{2}\right)  ^{2}}{c_{2}h_{0}\left(  1-c_{2}\right)
\left(  r+1\right)  }\right)
\end{align*}
If the latter expression is positive for $h_{0}>\bar{h}>0$, then $\lim
_{x_{0}\downarrow b_{+}}\Pi(x_{0},h_{0})$ is clearly negative. Therefore,
\[
\lim_{x_{0}\downarrow b_{+}}\Pi(x_{0},h_{0})<\max\left\{  0,\lim
_{x_{0}\downarrow b_{+}}\Pi(x_{0},\bar{h})\right\}  .
\]
But, using the definition $\bar{h}=\left(  c_{2}+r\right)  /\left(
1-c_{2}\right)  $ in (\ref{lim of Psi}), we obtain%
\[
\lim_{x_{0}\downarrow b_{+}}\Pi(x_{0},\bar{h})=-\frac{c_{1}\left(
1-c_{2}\right)  }{c_{2}\left(  1+r\right)  }-1+r\frac{\left(  r+c_{2}\right)
}{c_{2}\left(  r+1\right)  }=0.
\]
This implies that $\max_{\zeta\in\left[  \frac{h_{p}}{2}\tilde{x}_{0}%
,\frac{h_{p}}{2}\tilde{\lambda}_{1p}\right]  }\frac{\mathrm{d}}{\mathrm{d}%
\zeta}g_{p}\left(  \zeta\right)  $ is a.s. negative for sufficiently large $p$
and $\mathbf{n}$.

Now, since $\lambda_{p1}\overset{a.s.}{\rightarrow}x_{1}$, we have
\[
g_{p}^{^{\prime}}(0)\equiv\left.  -\frac{\mathrm{d}}{\mathrm{d}\zeta}%
g_{p}(\zeta)\right\vert _{\zeta=\frac{h_{p}}{2}\tilde{\lambda}_{p1}%
}\overset{a.s.}{\rightarrow}R_{0},
\]
where%
\[
R_{0}=\frac{c_{1}}{c_{2}\left(  h_{0}+1\right)  }+1+\frac{c_{1}}{h_{0}%
}(1+\alpha x_{1})+\frac{c_{1}}{\alpha h_{0}}(1+\alpha x_{1})^{2}m(x_{1}).
\]
Using (\ref{Iains result}), (\ref{def_x1}) and (\ref{def_H0}), we obtain%
\[
R_{0}=2c_{1}H_{0}/h_{0}.
\]

Exploiting the approach given in Olver (1997, pp. 81-82), we obtain%
\[
\int_{0}^{\frac{h_{p}}{2}\left(  \tilde{\lambda}_{p1}-\tilde{x}_{0}\right)
}e^{-n_{A}g_{p}(\zeta)}\psi\left(  \zeta\right)  x^{-\frac{1}{2}}%
\mathrm{d}x\overset{\mathrm{P}}{\sim}\left(  \frac{\pi}{n_{A}g_{p}^{^{\prime}%
}(0)}\right)  ^{\frac{1}{2}}e^{-n_{A}g_{p}\left(  \frac{h_{p}}{2}%
\tilde{\lambda}_{p1}\right)  }\psi\left(  \frac{h_{0}\alpha x_{1}}{2\left(
1+\alpha x_{1}\right)  }\right)  .
\]
On the other hand, direct calculation shows that
\[
\psi\left(  \frac{h_{0}\alpha x_{1}}{2\left(  1+\alpha x_{1}\right)  }\right)
=i\frac{\sqrt{2}\left(  c_{1}+c_{2}+c_{2}h_{0}\right)  }{\sqrt{c_{1}}%
\sqrt{c_{2}\mu_{0}^{2}+c_{1}^{2}-c_{1}+2c_{1}\mu_{0}}},
\]
and%
\[
g_{p}\left(  \frac{h_{p}}{2}\tilde{\lambda}_{p1}\right)  =\varphi\left(
\frac{h_{p}}{2}\tilde{\lambda}_{p1}\right)  +\frac{1}{2n_{A}}%
{\displaystyle\sum\nolimits_{j=2}^{p}}
\ln\left(  \tilde{\lambda}_{p1}-\tilde{\lambda}_{pj}\right)  .
\]
Therefore,%
\begin{equation}
J_{2,p}(\gamma,\Lambda)\overset{\mathrm{P}}{\sim}Z\left(  p,\mathbf{n}%
,h_{0}\right)  \left(  -\frac{\pi}{pH_{0}}\right)  ^{\frac{1}{2}}%
e^{-n_{A}\varphi\left(  \frac{h_{p}}{2}\tilde{\lambda}_{p1}\right)  }%
\prod_{j=2}^{p}\left(  \tilde{\lambda}_{p1}-\tilde{\lambda}_{pj}\right)
^{-\frac{1}{2}}, \label{J2}%
\end{equation}
where%
\[
Z\left(  p,\mathbf{n},h_{0}\right)  =\frac{C\left(  p,\mathbf{n}\right)
}{\sqrt{\pi p}}\frac{c_{1}+c_{2}\mu_{0}}{\sqrt{c_{2}\mu_{0}^{2}+c_{1}%
^{2}-c_{1}+2c_{1}\mu_{0}}}.
\]
As Lemma \ref{Lem:negl2} below shows, $J_{3,4,p}(\gamma,\Lambda)$ is
asymptotically dominated by $J_{2,p}(\theta,\lambda)$, which completes the
proof.$\square$

\begin{lemma}
\label{Lem:negl2}Under the hypothesis that $h_{p}=h_{0},$ uniformly in
$\gamma$ from any compact subset of $\mathbb{R}$%
\[
J_{3,4,p}(\gamma,\Lambda)=o_{\mathrm{P}}\left(  J_{2,p}(\gamma,\Lambda
)\right)  .
\]

\end{lemma}

\textbf{Proof:} Let us first consider the integral $J_{3,p}(\gamma,\Lambda)$
over the contour $\mathcal{C}_{3}$. For $z\in\mathcal{C}_{3}$, by definition,
we have $\mathfrak{R}\left(  \zeta\right)  \equiv\mathfrak{R}\left(
h_{p}z/2\right)  \geq-2u+v$. Therefore, the uniform approximation
(\ref{approx F11}) is still valid, and we have%
\begin{align}
J_{3,p}\left(  \gamma,\Lambda\right)  \overset{\mathrm{P}}{\sim}  &
\int_{\mathcal{C}_{3}}\frac{C\left(  p,\mathbf{n}\right)  e^{-n_{A}%
\varphi\left(  \zeta\right)  }}{\sqrt{2\pi n_{A}}i}\psi\left(  \zeta\right)
\prod_{j=1}^{p}\left(  z-\tilde{\lambda}_{pj}\right)  ^{-\frac{1}{2}%
}\mathrm{d}z\label{J3}\\
&  =\int_{\mathcal{C}_{3}}\frac{C\left(  p,\mathbf{n}\right)  e^{-n_{A}%
g_{p}\left(  \zeta\right)  }}{\sqrt{2\pi n_{A}}i}\psi\left(  \zeta\right)
\left(  z-\tilde{\lambda}_{p1}\right)  ^{-\frac{1}{2}}\mathrm{d}z.\nonumber
\end{align}
Let us show that, for $\zeta=h_{p}z/2$ with $z\in\mathcal{C}_{3},$%
\begin{equation}
\mathfrak{R}g_{p}\left(  \zeta\right)  >g_{p}\left(  h_{p}\tilde{x}%
_{0}/2\right)  . \label{minimum of real part}%
\end{equation}
Recall that%
\begin{align*}
g_{p}(\zeta)  &  =\varphi\left(  \zeta\right)  +\frac{1}{2n_{A}}%
{\displaystyle\sum\nolimits_{j=2}^{p}}
\ln\left(  \frac{2}{h_{p}}\zeta-\tilde{\lambda}_{pj}\right)  ,\text{ and}\\
\varphi\left(  \zeta\right)   &  =\left(  u-v\right)  \ln\left(  u-v\right)
+v\ln\left(  z_{+}+v\right)  -u\ln\left(  z_{+}+u\right)  -z_{+}.
\end{align*}
By definition of $\mathcal{C}_{3},$ as $z$ moves along $\mathcal{C}_{3}$ away
from $\tilde{x}_{0},$ $\zeta$ is changing so that $z_{+}$ moves along a circle
with center at $-u$ and radius $z_{0+}+u,$ where $z_{0+}$ is as defined in
(\ref{def z0+}). In particular, $\left\vert z_{+}+u\right\vert $ remains
constant, $\mathfrak{R}\left(  -z_{+}\right)  $ increases, and, since $v<u,$
$\left\vert z_{+}+v\right\vert $ increases too. Overall,%
\[
\mathfrak{R}\left(  \varphi\left(  \zeta\right)  \right)  =\left(  u-v\right)
\ln\left(  u-v\right)  +v\ln\left\vert z_{+}+v\right\vert -u\ln\left\vert
z_{+}+u\right\vert +\mathfrak{R}\left(  -z_{+}\right)
\]
is increasing. Note also that $\left\vert \zeta\right\vert =\left\vert
z_{+}\right\vert \left\vert z_{+}+v\right\vert /\left\vert z_{+}+u\right\vert
$ must increase, which implies that $\left\vert \frac{2}{h_{p}}\zeta
-\tilde{\lambda}_{pj}\right\vert $ is increasing for all $j\geq2,$ and thus,%
\[
\mathfrak{R}\left(  \frac{1}{2n_{A}}%
{\displaystyle\sum\nolimits_{j=2}^{p}}
\ln\left(  \frac{2}{h_{p}}\zeta-\tilde{\lambda}_{pj}\right)  \right)
\]
is increasing too. This implies (\ref{minimum of real part}).

On the other hand, in the above proof of Lemma \ref{Lem2a} we have shown that
$\frac{\mathrm{d}}{\mathrm{d}\zeta}g_{p}(\zeta)$ is continuous and
\textit{negative} on $\zeta\in\lbrack\frac{h_{p}}{2}\tilde{x}_{0},\frac{h_{p}%
}{2}\tilde{\lambda}_{p1}].$ Hence, there must exist $C>0$ such that, for any
$\zeta=h_{p}z/2$ with $z\in\mathcal{C}_{3},$
\[
\left\vert e^{-n_{A}g_{p}\left(  \zeta\right)  }\right\vert \leq e^{-n_{A}%
C}e^{-n_{A}\varphi\left(  \frac{h_{p}}{2}\tilde{\lambda}_{p1}\right)  }%
\prod_{j=2}^{p}\left(  \tilde{\lambda}_{p1}-\tilde{\lambda}_{pj}\right)
^{-\frac{1}{2}}.
\]
This inequality, together with (\ref{J2}) and (\ref{J3}) imply that%
\[
J_{3,p}\left(  \gamma,\Lambda\right)  =o_{\mathrm{P}}\left(  J_{2,p}%
(\gamma,\Lambda)\right)  .
\]
The fact that $J_{4,p}\left(  \gamma,\Lambda\right)  =o_{\mathrm{P}}\left(
J_{2,p}(\gamma,\Lambda)\right)  $ follows from pp. 29-31 of Johnstone and
Onatski (2014).$\square$

\section{Local Asymptotic Normality}

\subsection{Analysis for Setting 1}

Let us denote the likelihood ratio by
\begin{equation}
L_{p1}(\gamma,\Lambda)=\frac{f_{1}(\Lambda;h_{p})}{f_{1}(\Lambda;h_{0})}.
\label{like}%
\end{equation}
From Lemmas \ref{Lem1} and \ref{Lem2}, we obtain the following expression
\[
f_{1}(\Lambda;h_{p})=\frac{1}{2\pi i}C_{p1}(\Lambda)k_{p1}(h_{p})I_{p}%
(\gamma,\Lambda).
\]
Using Lemma \ref{Lem2}, we obtain%
\begin{equation}
L_{p1}(\gamma,\Lambda)\overset{\mathrm{P}}{\sim}\frac{k_{p1}(h_{p})}%
{k_{p1}(h_{0})}\left(  \frac{1-\frac{h_{p}}{1+h_{p}}\tilde{\lambda}_{p1}%
}{1-\frac{h_{0}}{1+h_{0}}\tilde{\lambda}_{p1}}\right)  ^{\frac{p-n-2}{2}}.
\label{LRsimp}%
\end{equation}
Consider a new local parameter%
\[
\theta_{1}=\gamma/\omega_{1}\left(  h_{0}\right)  ,
\]
where%
\[
\omega_{1}\left(  h_{0}\right)  =\frac{2h_{0}^{2}\left(  1+h_{0}\right)
^{2}r^{2}}{\left(  h_{0}-c_{2}\left(  1+h_{0}\right)  \right)  ^{2}}.
\]
We have the following lemma.

\begin{lemma}
\label{LAN 1}Let Under the null hypothesis that $h=h_{0}$, uniformly in
$\theta_{1}$ from any compact subset of $\mathbb{R}$,
\[
\ln L_{p1}(\gamma,\Lambda)=\theta_{1}\sqrt{p}\left(  \lambda_{p1}%
-x_{p1}\right)  -\frac{1}{2}\theta_{1}^{2}\tau_{1}^{2}(h_{0})+o_{\mathrm{P}%
}(1)
\]
where%
\begin{align*}
x_{p1}  &  =\frac{\left(  h_{0}+p/n_{1}\right)  \left(  h_{0}+1\right)
}{h_{0}-\left(  h_{0}+1\right)  p/n_{2}},\text{ and}\\
\tau_{1}^{2}(h_{0})  &  =2r^{2}\frac{h_{0}^{2}\left(  h_{0}+1\right)
^{2}\left(  h_{0}^{2}-c_{2}\left(  h_{0}+1\right)  ^{2}-c_{1}\right)
}{\left(  c_{2}-h_{0}+c_{2}h_{0}\right)  ^{4}}.
\end{align*}

\end{lemma}

\textbf{Proof:} Taking the logarithm of (\ref{LRsimp}) yields
\begin{align}
\ln L_{p1}(\gamma,\Lambda)  &  =\frac{n+2-p}{2}\left(  \ln\left(
1-\frac{\tilde{\lambda}_{p1}h_{0}}{1+h_{0}}\right)  -\ln\left(  1-\frac
{\tilde{\lambda}_{p1}h_{p}}{1+h_{p}}\right)  \right) \nonumber\\
&  -\left(  \frac{p-2}{2}\right)  \ln\frac{h_{p}}{h_{0}}+\left(  \frac
{p-n_{1}-2}{2}\right)  \ln\frac{1+h_{p}}{1+h_{0}}+o_{\mathrm{P}}(1).
\label{MLexp}%
\end{align}
Moreover, we have the following expansions
\begin{align}
\ln\left(  1-\frac{\tilde{\lambda}_{p1}h_{0}}{1+h_{0}}\right)  -\ln\left(
1-\frac{\tilde{\lambda}_{p1}h_{p}}{1+h_{p}}\right)   &  =p^{-\frac{1}{2}%
}\gamma\frac{\tilde{\lambda}_{p1}}{(1+h_{0})\left(  1+h_{0}(1-\tilde{\lambda
}_{p1})\right)  }\nonumber\\
&  -p^{-1}\gamma^{2}\frac{\tilde{\lambda}_{p1}}{(1+h_{0})^{2}(1+h_{0}%
(1-\tilde{\lambda}_{p1}))}\label{expa1}\\
&  \quad+p^{-1}\gamma^{2}\frac{{\tilde{\lambda}_{p1}}^{2}}{2(1+h_{0}%
)^{2}(1+h_{0}(1-\tilde{\lambda}_{p1}))^{2}}+o_{\mathrm{P}}(p^{-1}),\nonumber
\end{align}%
\begin{equation}
\ln\frac{1+h_{p}}{1+h_{0}}=p^{-\frac{1}{2}}\gamma\frac{1}{1+h_{0}}%
-p^{-1}\gamma^{2}\frac{1}{2(1+h_{0})^{2}}+o(p^{-1}), \label{expa2}%
\end{equation}
and
\begin{equation}
\ln\frac{h_{p}}{h_{0}}=p^{-\frac{1}{2}}\gamma h_{0}^{-1}-\frac{1}{2}%
p^{-1}\gamma^{2}h_{0}^{-2}+o(p^{-1}). \label{expa3}%
\end{equation}
Finally, using (\ref{expa1}), (\ref{expa2}), and (\ref{expa3}) in
(\ref{MLexp}) and noting the fact that $\lambda_{1}-x_{p1}%
\overset{a.s.}{\rightarrow}0,$ we obtain the statement of the lemma by
straightforward algebraic manipulations.$\square$

Lemma \ref{LAN 1} together with the asymptotic normality of $\sqrt{p}\left(
\lambda_{p1}-x_{p1}\right)  $ established in Proposition \ref{Proposition2}
imply, via Le Cam's First Lemma (see van der Vaart (1998), p.88), that the
sequences of the probability measures $\left\{  \mathbb{P}_{h_{0},p}\right\}
$ and $\left\{  \mathbb{P}_{h_{0}+\gamma/\sqrt{p},p}\right\}  $ describing the
joint distribution of the eigenvalues of $\mathbf{F}$ under the null
$H_{0}:h_{p}=h_{0}$ and under the local alternative $H_{1}:h_{p}=h_{0}%
+\gamma/\sqrt{p}$ are mutually contiguous. Moreover, the experiments $\left(
\mathbb{P}_{h_{0}+\theta_{1}\omega_{1}\left(  h_{0}\right)  /\sqrt{p}%
,p}:\theta_{1}\in\mathbb{R}\right)  $ converge to the Gaussian shift
experiment $\left(  N\left(  \theta_{1},\tau_{1}^{2}(h_{0})\right)
:\theta_{1}\in\mathbb{R}\right)  $. In particular, these experiments are
\textit{LAN}.

\subsection{Analysis for Setting 2}

Let us denote the likelihood ratio by%
\begin{equation}
L_{p2}(\gamma,\Lambda)=\frac{f_{2}(\Lambda;h_{p})}{f_{2}(\Lambda;h_{0})}.
\label{like2}%
\end{equation}
From Lemmas \ref{Lem1} and \ref{Lem2a}, we obtain the following expression
\[
f_{2}(\Lambda;h_{p})=\frac{1}{2\pi i}C_{p2}(\Lambda)k_{p2}(h_{p})J_{p}%
(\gamma,\Lambda).
\]
Using Lemma \ref{Lem2a}, and the definitions (\ref{def_fi}) and (\ref{def z+}%
), we obtain%
\begin{equation}
L_{p2}(\gamma,\Lambda)\overset{\mathrm{P}}{\sim}\exp\left[  -n_{A}%
{\displaystyle\sum\nolimits_{j=1}^{4}}
\left(  a_{j}\left(  h_{p}\right)  -a_{j}\left(  h_{0}\right)  \right)
\right]  , \label{LRsimp2}%
\end{equation}
where%
\[
a_{1}\left(  h\right)  =\frac{h+\ln h}{2},
\]%
\[
a_{2}(h)=-\frac{1}{2}\left(  \frac{h}{2}\tilde{\lambda}_{p1}-v+\sqrt{\left(
\frac{h}{2}\tilde{\lambda}_{p1}-v\right)  ^{2}+4u\frac{h}{2}\tilde{\lambda
}_{p1}}\right)  ,
\]%
\[
a_{3}(h)=-u\ln\left[  \frac{1}{2}\left(  \frac{h}{2}\tilde{\lambda}%
_{p1}-v+\sqrt{\left(  \frac{h}{2}\tilde{\lambda}_{p1}-v\right)  ^{2}%
+4u\frac{h}{2}\tilde{\lambda}_{p1}}\right)  \right]  ,
\]
and
\[
a_{4}\left(  h\right)  =\left(  u-v\right)  \ln\left[  \frac{1}{2}\left(
-\frac{h}{2}\tilde{\lambda}_{p1}-v+\sqrt{\left(  \frac{h}{2}\tilde{\lambda
}_{p1}-v\right)  ^{2}+4u\frac{h}{2}\tilde{\lambda}_{p1}}\right)  \right]  .
\]

We would like, first, to expand $a_{j}\left(  h_{p}\right)  -a_{j}(h_{0}),$
with $j=1,...,4,$ in the power series of $\gamma/\sqrt{p}$ up to, and
including, the terms of order $O_{\mathrm{P}}\left(  \frac{1}{p}\right)  .$
For $a_{1},$ we have%
\begin{equation}
a_{1}\left(  h_{p}\right)  -a_{1}(h_{0})=\frac{h_{0}+1}{2h_{0}}\frac{\gamma
}{\sqrt{p}}-\frac{1}{4h_{0}^{2}}\frac{\gamma^{2}}{p}. \label{a1}%
\end{equation}
For $a_{2},$ note that%
\begin{align*}
\left(  \frac{h_{p}}{2}\tilde{\lambda}_{p1}-v\right)  ^{2}+4u\frac{h_{p}}%
{2}\tilde{\lambda}_{p1}  &  =\left(  \frac{h_{0}}{2}\tilde{\lambda}%
_{p1}-v\right)  ^{2}+4u\frac{h_{0}}{2}\tilde{\lambda}_{p1}\\
&  +\left(  \frac{h_{0}}{2}\tilde{\lambda}_{p1}+2u-v\right)  \frac{\gamma
}{\sqrt{p}}\tilde{\lambda}_{p1}+\frac{\gamma^{2}}{4p}\tilde{\lambda}_{p1}^{2}.
\end{align*}
Using this expression and the facts that, when $h_{p}=h_{0},$
\[
\tilde{\lambda}_{p1}\overset{a.s.}{\rightarrow}\frac{\left(  h_{0}+1\right)
\left(  c_{1}+h_{0}\right)  }{h_{0}\left(  1+c_{1}/c_{2}+h_{0}\right)
},\text{ }u\rightarrow\frac{1+c_{1}/c_{2}-c_{1}}{2},\text{ and }%
v\rightarrow\frac{1-c_{1}}{2},
\]
we obtain after some algebra,%
\begin{equation}
a_{2}\left(  h_{p}\right)  -a_{2}(h_{0})=-\frac{\tilde{\lambda}_{1p}}%
{4}\left(  1+\frac{\frac{h_{0}}{2}\tilde{\lambda}_{p1}+2u-v}{S}\right)
\frac{\gamma}{\sqrt{p}}+C^{(2)}\frac{\gamma^{2}}{p}+o_{\mathrm{P}}\left(
\frac{1}{p}\right)  , \label{a2}%
\end{equation}
where%
\[
S=\sqrt{\left(  \frac{h_{0}}{2}\tilde{\lambda}_{p1}-v\right)  ^{2}%
+4u\frac{h_{0}}{2}\tilde{\lambda}_{p1}},
\]
and%
\[
C^{(2)}=\frac{c_{1}\left(  h_{0}+1\right)  ^{2}\left(  c_{1}+h_{0}\right)
^{2}\left(  c_{1}+c_{2}-c_{1}c_{2}\right)  \left(  c_{2}+c_{1}+h_{0}%
c_{2}\right)  }{2h_{0}^{2}\left(  2c_{2}h_{0}+c_{1}^{2}+c_{2}+c_{2}h_{0}%
^{2}+c_{1}+2h_{0}c_{1}\right)  ^{3}}.
\]

For $a_{3},$ we have%
\begin{align}
a_{3}\left(  h_{p}\right)  -a_{3}(h_{0})  &  =-\frac{u\tilde{\lambda}%
_{p1}\left(  \frac{h_{0}}{2}\tilde{\lambda}_{p1}+2u-v+S\right)  }{2\left(
\left(  \frac{h_{0}}{2}\tilde{\lambda}_{p1}-v\right)  ^{2}+4u\frac{h_{0}}%
{2}\tilde{\lambda}_{p1}+\left(  \frac{h_{0}}{2}\tilde{\lambda}_{p1}-v\right)
S\right)  }\frac{\gamma}{\sqrt{p}}\label{a3}\\
&  +C^{(2)}C^{(3)}\frac{\gamma^{2}}{p}+o_{\mathrm{P}}\left(  \frac{1}%
{p}\right)  ,\nonumber
\end{align}
where%
\[
C^{(3)}=\frac{\left(  c_{1}+c_{2}-c_{1}c_{2}\right)  }{c_{2}\left(
c_{1}+h_{0}\right)  }+\frac{\left(  2c_{2}h_{0}+c_{1}^{2}+c_{2}+c_{2}h_{0}%
^{2}+c_{1}+2h_{0}c_{1}\right)  \left(  c_{1}+c_{2}+c_{2}h_{0}\right)  }%
{2c_{1}c_{2}\left(  c_{1}+h_{0}\right)  ^{2}}.
\]
Finally, for $a_{4},$ we obtain%
\begin{align}
a_{4}\left(  h_{p}\right)  -a_{4}(h_{0})  &  =\frac{\left(  u-v\right)
\tilde{\lambda}_{p1}\left(  \frac{h_{0}}{2}\tilde{\lambda}_{p1}+2u-v-S\right)
}{2\left(  \left(  \frac{h_{0}}{2}\tilde{\lambda}_{p1}-v\right)  ^{2}%
+4u\frac{h_{0}}{2}\tilde{\lambda}_{p1}+\left(  -\frac{h_{0}}{2}\tilde{\lambda
}_{p1}-v\right)  S\right)  }\frac{\gamma}{\sqrt{p}}\label{a4}\\
&  -C^{(2)}C^{(4)}\frac{\gamma^{2}}{p}+o_{\mathrm{P}}\left(  \frac{1}%
{p}\right)  ,\nonumber
\end{align}
where%
\[
C^{(4)}=\frac{c_{1}+c_{2}+c_{2}h_{0}}{c_{2}\left(  c_{1}+h_{0}\right)  }%
+\frac{\left(  c_{1}+c_{2}-c_{1}c_{2}\right)  \left(  c_{1}+c_{2}+c_{2}%
h_{0}^{2}+c_{1}^{2}+2c_{1}h_{0}+2c_{2}h_{0}\right)  }{2c_{2}\left(
c_{1}+h_{0}\right)  ^{2}\left(  c_{1}+c_{2}+c_{2}h_{0}\right)  }.
\]
Summing up the $\gamma^{2}/p$ terms in the expansions (\ref{a1}-\ref{a4}), we
obtain that the $\gamma^{2}/p$ term in the expansion of $%
{\displaystyle\sum\nolimits_{j=1}^{4}}
\left(  a_{j}\left(  h_{p}\right)  -a_{j}\left(  h_{0}\right)  \right)  $,
which we will refer as $T_{2}$, equals
\begin{equation}
T_{2}=-\frac{1}{4}c_{1}\frac{c_{1}+c_{2}+c_{2}h_{0}^{2}-h_{0}^{2}+2c_{2}h_{0}%
}{h_{0}^{2}\left(  c_{1}+c_{2}+c_{2}h_{0}^{2}+c_{1}^{2}+2c_{1}h_{0}%
+2c_{2}h_{0}\right)  }\frac{\gamma^{2}}{p}. \label{T2}%
\end{equation}

Now let $\Delta=\sqrt{p}\left(  \lambda_{p1}-x_{p1}\right)  ,$ where%
\[
x_{p1}=\frac{\left(  h_{0}+p/n_{1}\right)  \left(  h_{0}+1\right)  }%
{h_{0}-\left(  h_{0}+1\right)  p/n_{2}},
\]
so that%
\begin{align*}
\tilde{\lambda}_{p1}  &  =\frac{\lambda_{p1}}{\frac{n_{2}}{n_{1}}+\lambda
_{p1}}=\frac{x_{p1}+\Delta/\sqrt{p}}{\frac{n_{2}}{n_{1}}+x_{p1}+\Delta
/\sqrt{p}}\\
&  =\frac{\left(  h_{0}+1\right)  \left(  p/n_{1}+h_{0}\right)  }{h_{0}\left(
1+n_{2}/n_{1}+h_{0}\right)  }+\frac{\Delta}{\sqrt{p}}\frac{c_{2}c_{1}\left(
c_{2}+c_{2}h_{0}-h_{0}\right)  ^{2}}{h_{0}^{2}\left(  c_{2}+c_{1}+c_{2}%
h_{0}\right)  ^{2}}+o_{\mathrm{P}}\left(  \frac{1}{\sqrt{p}}\right)  ,
\end{align*}
Our next goal is to expand the weights on $\gamma/\sqrt{p}$ in expansions
(\ref{a1}-\ref{a4}) into power series of $\Delta/\sqrt{p}$ up to the linear
term only.

For (\ref{a2}), we have%
\[
-\frac{\tilde{\lambda}_{p1}}{4}\left(  1+\frac{\frac{h_{0}}{2}\tilde{\lambda
}_{p1}+2u-v}{S}\right)  =\tau_{0}^{(2)}+\left[  \tau_{11}^{(2)}+\tau
_{12}^{(2)}\right]  \frac{\Delta}{\sqrt{p}}+o_{\mathrm{P}}\left(  \frac
{1}{\sqrt{p}}\right)  ,
\]
where%

\[
\tau_{11}^{(2)}=-\frac{c_{1}\left(  c_{2}+c_{2}h_{0}-h_{0}\right)  ^{2}%
}{2h_{0}^{2}\left(  c_{1}+c_{2}+c_{2}h_{0}^{2}+c_{1}^{2}+2c_{1}h_{0}%
+2c_{2}h_{0}\right)  },
\]%
\[
\tau_{12}^{(2)}=\frac{c_{1}^{2}\left(  h_{0}+1\right)  \left(  c_{1}%
+h_{0}\right)  \left(  c_{1}+c_{2}-c_{1}c_{2}\right)  \left(  c_{2}+c_{2}%
h_{0}-h_{0}\right)  ^{2}}{h_{0}^{2}\left(  c_{1}+c_{2}+c_{2}h_{0}^{2}%
+c_{1}^{2}+2c_{1}h_{0}+2c_{2}h_{0}\right)  ^{3}},
\]
and $\tau_{0}^{(2)}$ is a complicated function of $h_{0},p,n_{1},$ and $n_{2}%
$, which we do not report here.

For (\ref{a3}), we have
\[
-\frac{u\tilde{\lambda}_{p1}\left(  \frac{h_{0}}{2}\tilde{\lambda}%
_{p1}+2u-v+S\right)  }{2\left(  \left(  \frac{h_{0}}{2}\tilde{\lambda}%
_{p1}-v\right)  ^{2}+4u\frac{h_{0}}{2}\tilde{\lambda}_{p1}+\left(  \frac
{h_{0}}{2}\tilde{\lambda}_{p1}-v\right)  S\right)  }=\tau_{0}^{(3)}+\left[
\tau_{11}^{(3)}+\tau_{12}^{(3)}\right]  \frac{\Delta}{\sqrt{p}}+o_{\mathrm{P}%
}\left(  \frac{1}{\sqrt{p}}\right)  ,
\]
where%
\[
\tau_{11}^{(3)}=\frac{c_{1}\left(  c_{2}-h_{0}+c_{2}h_{0}\right)  ^{2}\left(
c_{1}+c_{2}-c_{1}c_{2}\right)  \left(  c_{1}+c_{2}+c_{2}h_{0}\right)  \left(
c_{1}+c_{2}+c_{2}h_{0}^{2}-c_{1}^{2}+2c_{2}h_{0}\right)  }{2h_{0}^{2}%
c_{2}\left(  c_{1}+c_{2}+c_{2}h_{0}^{2}+c_{1}^{2}+2c_{1}h_{0}+2c_{2}%
h_{0}\right)  ^{3}},
\]%
\[
\tau_{12}^{(3)}=-\frac{c_{1}\left(  h_{0}+1\right)  \left(  c_{2}-h_{0}%
+c_{2}h_{0}\right)  ^{2}\left(  c_{1}+c_{2}-c_{1}c_{2}\right)  }{2h_{0}%
^{2}\left(  c_{1}+c_{2}+c_{2}h_{0}^{2}+c_{1}^{2}+2c_{1}h_{0}+2c_{2}%
h_{0}\right)  ^{2}},
\]
and $\tau_{0}^{(3)}$ is a complicated function of $h_{0},p,n_{1},$ and $n_{2}%
$, which we do not report here.

For (\ref{a4}), we have%
\[
\frac{\left(  u-v\right)  \tilde{\lambda}_{p1}\left(  \frac{h_{0}}{2}%
\tilde{\lambda}_{p1}+2u-v-S\right)  }{2\left(  \left(  \frac{h_{0}}{2}%
\tilde{\lambda}_{p1}-v\right)  ^{2}+4u\frac{h_{0}}{2}\tilde{\lambda}%
_{p1}+\left(  -\frac{h_{0}}{2}\tilde{\lambda}_{p1}-v\right)  S\right)  }%
=\tau_{0}^{(4)}+\left[  \tau_{11}^{(4)}+\tau_{12}^{(4)}\right]  \frac{\Delta
}{\sqrt{p}}+o_{\mathrm{P}}\left(  \frac{1}{\sqrt{p}}\right)  ,
\]
where%
\[
\tau_{11}^{(4)}=\frac{c_{1}^{3}\left(  c_{2}-h_{0}+c_{2}h_{0}\right)
^{2}\left(  c_{1}+c_{2}-c_{1}c_{2}\right)  \left(  -c_{1}-c_{2}+c_{2}h_{0}%
^{2}+c_{1}^{2}+2c_{1}c_{2}+2c_{1}c_{2}h_{0}\right)  }{2c_{2}h_{0}^{2}\left(
c_{1}+c_{2}+c_{2}h_{0}\right)  \left(  c_{1}+c_{2}+c_{2}h_{0}^{2}+c_{1}%
^{2}+2c_{1}h_{0}+2c_{2}h_{0}\right)  ^{3}},
\]%
\[
\tau_{12}^{(4)}=-\frac{c_{1}^{2}\left(  h_{0}+1\right)  \left(  c_{2}%
-h_{0}+c_{2}h_{0}\right)  ^{2}\left(  c_{1}+c_{2}-c_{1}c_{2}\right)  }%
{2h_{0}^{2}\left(  c_{1}+c_{2}+c_{2}h_{0}\right)  \left(  c_{1}+c_{2}%
+c_{2}h_{0}^{2}+c_{1}^{2}+2c_{1}h_{0}+2c_{2}h_{0}\right)  ^{2}},
\]
and $\tau_{0}^{(4)}$ is a complicated function of $h_{0},p,n_{1},$ and $n_{2}%
$, which we do not report here.

We have verified, using Maple symbolic algebra software, that%
\[
\tau^{(2)}+\tau^{(3)}+\tau^{(4)}=-\frac{1}{2}\frac{1+h_{0}}{h_{0}},
\]
which is exactly the negative of the term on $\gamma/\sqrt{p}$ in (\ref{a1}).
Hence, the term on $\gamma/\sqrt{p}$ in the expansion of $%
{\displaystyle\sum\nolimits_{j=1}^{4}}
\left(  a_{j}\left(  h_{p}\right)  -a_{j}\left(  h_{0}\right)  \right)  $ is
zero. Further, we have verified that%
\[%
{\displaystyle\sum\nolimits_{j=2}^{4}}
\left(  \tau_{11}^{(j)}+\tau_{12}^{(j)}\right)  =-\frac{1}{2}c_{1}%
\frac{\left(  c_{2}-h_{0}+c_{2}h_{0}\right)  ^{2}}{h_{0}^{2}\left(
c_{1}+c_{2}+c_{2}h_{0}^{2}+c_{1}^{2}+2c_{1}h_{0}+2c_{2}h_{0}\right)  }.
\]
This equality, together with (\ref{LRsimp2}) and (\ref{T2}) imply that%
\begin{align}
&  \ln L_{p1}(\gamma,\Lambda)\overset{\mathrm{P}}{\sim}\frac{1}{2}%
\frac{\left(  c_{2}-h_{0}+c_{2}h_{0}\right)  ^{2}}{h_{0}^{2}\left(
c_{1}+c_{2}+c_{2}h_{0}^{2}+c_{1}^{2}+2c_{1}h_{0}+2c_{2}h_{0}\right)  }%
\gamma\Delta\label{expansion}\\
&  +\frac{1}{4}\frac{c_{1}+c_{2}+c_{2}h_{0}^{2}-h_{0}^{2}+2c_{2}h_{0}}%
{h_{0}^{2}\left(  c_{1}+c_{2}+c_{2}h_{0}^{2}+c_{1}^{2}+2c_{1}h_{0}+2c_{2}%
h_{0}\right)  }\gamma^{2}.\nonumber
\end{align}

Consider a different local parameter%
\[
\theta_{2}=\gamma/\omega_{2}\left(  h_{0}\right)  ,
\]
where%
\[
\omega_{2}\left(  h_{0}\right)  =\frac{2h_{0}^{2}\left(  c_{1}+c_{2}%
+c_{2}h_{0}^{2}+c_{1}^{2}+2c_{1}h_{0}+2c_{2}h_{0}\right)  }{\left(
h_{0}-c_{2}\left(  1+h_{0}\right)  \right)  ^{2}}.
\]
Asymptotic approximation (\ref{expansion}) implies the following lemma.

\begin{lemma}
\label{LAN 2}Under the null hypothesis that $h=h_{0}$, uniformly in
$\theta_{2}$ from any compact subset of $\mathbb{R}$,
\[
\ln L_{p2}(\gamma,\Lambda)=\theta_{2}\sqrt{p}\left(  \lambda_{p1}%
-x_{p1}\right)  -\frac{1}{2}\theta_{2}^{2}\tau_{2}^{2}(h_{0})+o_{\mathrm{P}%
}(1)
\]
where%
\begin{align*}
x_{p1}  &  =\frac{\left(  h_{0}+p/n_{1}\right)  \left(  h_{0}+1\right)
}{h_{0}-\left(  h_{0}+1\right)  p/n_{2}},\text{ and}\\
\tau_{2}^{2}(h_{0})  &  =\frac{2h_{0}^{2}\left(  h_{0}^{2}-c_{2}\left(
1+h_{0}\right)  ^{2}-c_{1}\right)  \left(  \left(  c_{1}+c_{2}\right)  \left(
1+h_{0}\right)  ^{2}-c_{1}\left(  h_{0}^{2}-c_{1}\right)  \right)  }{\left(
c_{2}-h_{0}+c_{2}h_{0}\right)  ^{4}}.
\end{align*}

\end{lemma}

Similarly to the case of Setting 1, Lemma \ref{LAN 2} together with the
asymptotic normality of $\sqrt{p}\left(  \lambda_{p1}-x_{p1}\right)  $
established in Proposition \ref{Proposition2} imply, via Le Cam's First Lemma
(see van der Vaart (1998), p.88), that the sequences of the probability
measures $\left\{  \mathbb{P}_{h_{0},p}\right\}  $ and $\left\{
\mathbb{P}_{h_{0}+\gamma/\sqrt{p},p}\right\}  $ describing the joint
distribution of the eigenvalues of $\mathbf{F}$ under the null $H_{0}%
:h_{p}=h_{0}$ and under the local alternative $H_{1}:h_{p}=h_{0}+\gamma
/\sqrt{p}$ are mutually contiguous. Moreover, the experiments $\left(
\mathbb{P}_{h_{0}+\theta_{2}\omega_{2}\left(  h_{0}\right)  /\sqrt{p}%
,p}:\theta_{2}\in\mathbb{R}\right)  $ converge to the Gaussian shift
experiment $\left(  N\left(  \theta_{2},\tau_{2}^{2}(h_{0})\right)
:\theta_{2}\in\mathbb{R}\right)  $. In particular, these experiments are
\textit{LAN}.

\section{Conclusion}

In this paper, we establish the Local Asymptotic Normality of the experiments
of observing the eigenvalues of the F-ratio $\mathbf{F}\equiv\left(
B/n_{2}\right)  ^{-1}A/n_{A}$ of two large-dimensional Wishart matrices. The
experiments are parameterized by the value of a single spike that describes
the \textquotedblleft ratio\textquotedblright\ of the covariance parameters of
$A$ and $B$, or, in the case of equal covariance parameters, the
non-centrality parameter of $A$. We find that the asymptotic behavior of the
log ratio of the joint density of the eigenvalues of $\mathbf{F},$ which
corresponds to a super-critical spike, to their joint density under a local
deviation from this value depends only on the largest eigenvalue $\lambda
_{p1}$. This implies, in particular, that the best statistical inference about
a super-critical spike in the local asymptotic regime is based on the largest
eigenvalue only.

As a by-product of our analysis, in a multi-spike setting, we establish the
joint asymptotic normality of a few of the largest eigenvalues of $\mathbf{F}$
that correspond to the super-critical spikes. We derive an explicit formulas
for the almost sure limits of these eigenvalues, and for the asymptotic
variances of their fluctuations around these limits.

\section{Acknowledgements}

\label{sec:acknowledgements}

This work was supported in part by NIH grant 5R01 EB 001988 (PD, IMJ), the
Simons Foundation Math + X program (PD), NSF grant DMS 1407813 (IMJ), and the
the J.M. Keynes Fellowships Fund, University of Cambridge (AO).

\section{Appendix}

\subsection{Proof of Lemma \ref{Lemma 13}}

We will need the following two lemmas.

\begin{lemma}
\label{McLeish}(McLeish 1974) Let $\left\{  X_{pj},\mathcal{G}_{pj}%
,j=1,...,p\right\}  $ be a martingale difference array on the probability
triple $\left(  \Omega,\mathcal{G},P\right)  $. If the following conditions
are satisfied: a) Lindeberg's condition: for all $\varepsilon>0$, $%
{\displaystyle\sum\nolimits_{j}}
\int_{\left\vert X_{pj}\right\vert >\varepsilon}X_{pj}^{2}\mathrm{d}%
P\rightarrow0$ as $p\rightarrow\infty$; b) $%
{\displaystyle\sum\nolimits_{j}}
X_{pj}^{2}\overset{\mathrm{P}}{\rightarrow}1,$ then $%
{\displaystyle\sum\nolimits_{j}}
X_{pj}\overset{d}{\rightarrow}N\left(  0,1\right)  $.
\end{lemma}

\textbf{Proof:} This is a consequence of Theorem (2.3) of McLeish (1974). Two
conditions of the theorem: i) $\max_{j\leq p}\left\vert X_{pj}\right\vert $ is
uniformly bounded in $L_{2}$ norm, and ii) $\max_{j\leq p}\left\vert
X_{pj}\right\vert \overset{\mathrm{P}}{\rightarrow}0$, are replaced here by
the Lindeberg condition.$\square$

\begin{lemma}
\label{Hall and Heyde}(Hall and Heyde) Let $\left\{  X_{pj},\mathcal{G}%
_{pj},j=1,...,p\right\}  $ be a martingale difference array, and define
$V_{pJ}^{2}=%
{\displaystyle\sum\nolimits_{j=1}^{J}}
E\left(  X_{pj}^{2}|\mathcal{G}_{p,j-1}\right)  $ and $U_{pJ}^{2}=%
{\displaystyle\sum\nolimits_{j=1}^{J}}
X_{pj}^{2}$ for $J=1,...,p$. Suppose that the conditional variances
$V_{pp}^{2}$ are tight, that is $\sup_{p}P\left(  V_{pp}^{2}>\varepsilon
\right)  \rightarrow0$ as $\varepsilon\rightarrow\infty$, and that the
conditional Lindeberg condition holds, that is, for all $\varepsilon>0$, $%
{\displaystyle\sum\nolimits_{j}}
E\left[  X_{pj}^{2}\mathbf{1}\left\{  \left\vert X_{pj}\right\vert
>\varepsilon\right\}  |\mathcal{G}_{p,j-1}\right]  \overset{\mathrm{P}%
}{\rightarrow}0$. Then $\max_{J}\left\vert U_{pJ}^{2}-V_{pJ}^{2}\right\vert
\overset{\mathrm{P}}{\rightarrow}0.$
\end{lemma}

\textbf{Proof:} This is a shortened version of Theorem 2.23 in Hall and Heyde
(1980).$\square$

Let $f_{q}\left(  \lambda\right)  ,$ $q=1,...,Q,$ be such that $f_{q}%
(\lambda)=g_{q}\left(  \lambda\right)  $ for $\lambda\in\left[  l_{i}%
,L_{i}\right]  $ and $f_{q}(\lambda)=0$ otherwise. Consider random variables
\[
X_{pj}=\frac{1}{\sqrt{p}}%
{\displaystyle\sum\nolimits_{\left(  q,s,t\right)  \in\Theta}}
\gamma_{qst}f_{q}\left(  \lambda_{pj}^{(i)}\right)  \left(  \zeta_{js}%
\zeta_{jt}-\delta_{st}\right)  ,
\]
where $\gamma_{qst}$ are some constants. Let $\mathcal{G}_{pJ}$ be the
$\sigma$-algebra generated by $\lambda_{p1}^{(i)},...,\lambda_{pp}^{(i)}$ and
$\zeta_{js}$ with $j=1,...,J;$ $s=1,...,m$. Clearly, $\left\{  X_{pj}%
,\mathcal{G}_{pj},j=1,...,p\right\}  $ form a martingale difference array. Let
$K$ be the number of different triples $\left(  q,s,t\right)  \in\Theta.$
Consider an arbitrary order in $\Theta$. In H\"{o}lder's inequality%
\[%
{\displaystyle\sum\nolimits_{a=1}^{K}}
y_{a}z_{a}\leq\left(
{\displaystyle\sum\nolimits_{a=1}^{K}}
\left(  y_{a}\right)  ^{b}\right)  ^{1/b}\left(
{\displaystyle\sum\nolimits_{a=1}^{K}}
\left(  z_{a}\right)  ^{c}\right)  ^{1/c},
\]
which holds for $y_{a}>0,$ $z_{a}>0$, $b>1,$ $c>1,$ and $1/b+1/c=1,$ take%
\[
y_{a}=\left\vert \frac{1}{\sqrt{p}}\gamma_{qst}f_{q}\left(  \lambda_{pj}%
^{(i)}\right)  \left(  \zeta_{js}\zeta_{jt}-\delta_{st}\right)  \right\vert ,
\]
where $\left(  q,s,t\right)  $ is the $a$-th triple in $\Theta,$ $z_{a}=1,$
and $b=2+\delta$ for some $\delta>0$. Then, the inequality implies that%
\begin{equation}
\left\vert X_{pj}\right\vert ^{2+\delta}\leq K^{1+\delta}R_{i}^{2+\delta}%
{\displaystyle\sum\nolimits_{\left(  q,s,t\right)  \in\Theta}}
\left\vert \frac{1}{\sqrt{p}}\gamma_{qst}\left(  \zeta_{js}\zeta_{jt}%
-\delta_{st}\right)  \right\vert ^{2+\delta}, \label{Holder}%
\end{equation}
where%
\[
R_{i}=\max_{q=1,...,Q}\sup_{\lambda\in\left[  l_{i},L_{i}\right]  }\left\vert
g_{q}\left(  \lambda\right)  \right\vert .
\]
Since $\zeta_{js}$ are i.i.d. $N(0,1),$ (\ref{Holder}) implies that $%
{\displaystyle\sum\nolimits_{j=1}^{p}}
E\left\vert X_{pj}\right\vert ^{2+\delta}\rightarrow0$ as $p\rightarrow
\infty,$ which means that the Lyapunov condition holds for $X_{pj}$. As is
well known, Lyapunov's condition implies Lindeberg's condition. Hence,
condition a) of Lemma \ref{McLeish} is satisfied for $X_{pj}$.

Let us consider $%
{\displaystyle\sum\nolimits_{j=1}^{p}}
X_{pj}^{2}$. Since the convergence in mean implies the convergence in
probability, the conditional Lindeberg condition is satisfied for $X_{pj}$
because the unconditional Lindeberg condition is satisfied as checked above.
Further, in notations of Lemma \ref{Hall and Heyde}, it is easy to see that
\[
V_{pp}^{2}=%
{\displaystyle\sum\nolimits_{q,q_{1}}}
\left[  \left(
{\displaystyle\sum\nolimits_{1\leq s\leq t\leq m}}
\gamma_{qst}\gamma_{q_{1}st}\left(  1+\delta_{st}\right)  \right)  \frac{1}{p}%
{\displaystyle\sum\nolimits_{j=1}^{p}}
f_{q}\left(  \lambda_{pj}^{(i)}\right)  f_{q_{1}}\left(  \lambda_{pj}%
^{(i)}\right)  \right]  .
\]
The convergence of the empirical distribution of $\lambda_{p1}^{(i)}%
,...,\lambda_{pp}^{(i)}$ to $G_{x_{i}}$ and the equality of $g_{q}$ and
$f_{q}$ on the support of $G_{x_{i}}$ implies that%
\[
V_{pp}^{2}\overset{\mathrm{P}}{\rightarrow}\Sigma\equiv%
{\displaystyle\sum\nolimits_{q,q_{1}}}
\left[  \left(
{\displaystyle\sum\nolimits_{1\leq s\leq t\leq m}}
\gamma_{qst}\gamma_{q_{1}st}\left(  1+\delta_{st}\right)  \right)  \int
g_{q}\left(  \lambda\right)  g_{q_{1}}\left(  \lambda\right)  \mathrm{d}%
G_{x_{i}}\right]  .
\]
In particular, $V_{pp}^{2}$ is tight and Lemma \ref{Hall and Heyde} applies.
Therefore, $%
{\displaystyle\sum\nolimits_{j=1}^{p}}
X_{pj}^{2}$ converges to the same limit as $V_{pp}^{2}$. Thus, by Lemma
\ref{McLeish}, we get $%
{\displaystyle\sum\nolimits_{j=1}^{p}}
X_{pj}\overset{d}{\rightarrow}N(0,\Sigma).$

Finally, let
\[
Y_{pj}=\frac{1}{\sqrt{p}}%
{\displaystyle\sum\nolimits_{\left(  q,s,t\right)  \in\Theta}}
\gamma_{qst}g_{q}\left(  \lambda_{pj}^{(i)}\right)  \left(  \zeta_{js}%
\zeta_{jt}-\delta_{st}\right)  .
\]
Since
\[
\Pr\left(
{\displaystyle\sum\nolimits_{j=1}^{p}}
X_{pj}\neq%
{\displaystyle\sum\nolimits_{j=1}^{p}}
Y_{pj}\right)  \rightarrow0
\]
as $p\rightarrow\infty$, we have $%
{\displaystyle\sum\nolimits_{j=1}^{p}}
Y_{pj}\overset{d}{\rightarrow}N(0,\Sigma)$. Lemma \ref{Lemma 13} follows from
this convergence via the Cramer-Wold device.$\square$

\subsection{Derivation of (\ref{ingredient1}), (\ref{ingredient2}), and
(\ref{ingredient3})}

Expression (\ref{ingredient1}) immediately follows from (\ref{root}). Next,
differentiating identity (\ref{eq3}) with respect to $z$, we obtain%
\[
1+\frac{c_{1}m_{x}^{\prime}\left(  z\right)  }{\left(  1+c_{1}m_{x}\left(
z\right)  \right)  ^{2}}=\frac{m_{x}^{\prime}\left(  z\right)  }{m_{x}%
^{2}\left(  z\right)  }+\frac{-x^{2}c_{2}m_{x}^{\prime}\left(  z\right)
}{\left(  1-c_{2}xm_{x}\left(  z\right)  \right)  ^{2}}.
\]
Setting $z=0$ and $x=x_{i},$ and using the fact that
\begin{equation}
m_{x_{i}}\left(  0\right)  =-\left(  h_{i}+c_{1}\right)  ^{-1}, \label{mx0}%
\end{equation}
which follows from (\ref{root}), we obtain%
\[
1+\frac{c_{1}m_{x_{i}}^{\prime}\left(  0\right)  }{\left(  1-c_{1}\left(
h_{i}+c_{1}\right)  ^{-1}\right)  ^{2}}=\frac{m_{x_{i}}^{\prime}\left(
0\right)  }{\left(  h_{i}+c_{1}\right)  ^{-2}}+\frac{-x_{i}^{2}c_{2}m_{x_{i}%
}^{\prime}\left(  0\right)  }{\left(  1+c_{2}x_{i}\left(  h_{i}+c_{1}\right)
^{-1}\right)  ^{2}}.
\]
Using the definition (\ref{x0}) of $x_{i}$, we obtain%
\begin{align*}
1+\frac{c_{1}m_{x_{i}}^{\prime}\left(  0\right)  }{\left(  1-c_{1}\left(
h_{i}+c_{1}\right)  ^{-1}\right)  ^{2}}  &  =\frac{m_{x_{i}}^{\prime}\left(
0\right)  }{\left(  h_{i}+c_{1}\right)  ^{-2}}\\
&  -\frac{\left(  h_{i}+c_{1}\right)  ^{2}\left(  h_{i}+1\right)  ^{2}%
c_{2}m_{x_{0}}^{\prime}\left(  0\right)  }{h_{i}^{2}},
\end{align*}
which implies (\ref{ingredient2}). Finally, differentiating identity
(\ref{eq3}) with respect to $x$, we obtain%
\begin{align*}
\frac{c_{1}\mathrm{d}m_{x}\left(  z\right)  /\mathrm{d}x}{\left(  1+c_{1}%
m_{x}\left(  z\right)  \right)  ^{2}}  &  =\frac{\mathrm{d}m_{x}\left(
z\right)  /\mathrm{d}x}{\left(  m_{x}\left(  z\right)  \right)  ^{2}}\\
&  +\frac{-1+c_{2}xm_{x}\left(  z\right)  -x\left(  c_{2}m_{x}\left(
z\right)  +c_{2}x\mathrm{d}m_{x}\left(  z\right)  /\mathrm{d}x\right)
}{\left(  1-c_{2}xm_{x}\left(  z\right)  \right)  ^{2}}.
\end{align*}
Setting $z=0$ and $x=x_{i},$ we obtain%
\[
\frac{c_{1}\mathrm{d}m_{x_{i}}\left(  0\right)  /\mathrm{d}x}{\left(
1+c_{1}m_{x_{i}}\left(  0\right)  \right)  ^{2}}=\frac{\mathrm{d}m_{x_{i}%
}\left(  0\right)  /\mathrm{d}x}{\left(  m_{x_{i}}\left(  0\right)  \right)
^{2}}+\frac{-1-c_{2}x_{i}^{2}\mathrm{d}m_{x_{i}}\left(  0\right)
/\mathrm{d}x}{\left(  1-c_{2}x_{i}m_{x_{i}}\left(  0\right)  \right)  ^{2}}.
\]
This equality, the definition (\ref{x0}) of $x_{i},$ and equation (\ref{mx0})
imply (\ref{ingredient3}).

\end{document}